\documentclass{amsart}

\usepackage{amsthm,amsfonts,amsmath,amssymb,latexsym,epsfig}
\usepackage{upref,amssymb,eucal,ae}
\usepackage[all,cmtip]{xy}

\newtheorem{theorem}{Theorem}
\newtheorem{lemma}[theorem]{Lemma}

\newtheorem{corollary}[theorem]{Corollary}
\newenvironment{example}{\medskip \refstepcounter{theorem}
\noindent  {\bf Example \thetheorem}.\rm}{\,}

\makeatletter
\def\Ddots{\mathinner{\mkern1mu\raise\p@
\vbox{\kern7\p@\hbox{.}}\mkern2mu
\raise4\p@\hbox{.}\mkern2mu\raise7\p@\hbox{.}\mkern1mu}}
\makeatother

\newcounter{spthe}

\def\<{\langle}
\def\a{\alpha}
\def\>{\rangle}

\def\tm{\tilde{M}}
\def\tn{\tilde{n}}
\def\tg{\tilde{g}}
\def\mb#1{{\mathbb #1}}
\def\mc#1{{\mathcal #1}}

\def\kuno{\bigcirc \hspace{-3.75mm} \wedge \hspace{1.4mm}}
 
\begin{document}
\title[$\mb{S}^6$ is not diffeomorphic to a complex manifold]{$\mb{S}^6$ (or 
any of $\mb{S}^2 \times \mb{S}^4$, $\mb{S}^2\times\mb{S}^6$, or 
$\mb{S}^6\times \mb{S}^6$, respectively) is not  diffeomorphic to a complex 
manifold}
\author{Santiago R. Simanca}
\thanks{This work was carried out, in great part, while the author was 
supported by the Simons Foundation Visiting Professorship award 657746 at CIMS}
\email{srsimanca@gmail.com}

\begin{abstract}
We identify all metrics on a closed $n$-manifold with their Nash isometric 
embeddings into a standard sphere of large, but fixed dimension, and use 
the Palais' isotopic extension theorem to identify their deformations 
with the isotopic deformations of their embeddings, the deformations of 
metrics in a conformal class identified with their corresponding isotopic 
conformal deformations. If $n\geq 3$, we characterize metrics of 
constant scalar curvature in terms of properties of extrinsic quantities of 
their associated embeddings, and prove that any metric on the manifold of 
constant positive scalar curvature, which can be minimally embedded 
into this background sphere, is a Yamabe metric in its conformal class.
We then use Simons' gap theorem to study the extrinsic quantities of 
almost complex Hermitian deformations, by Yamabe metrics, of the standard 
minimal almost complex isometric embeddings of $\mb{S}^6$, 
$\mb{S}^2 \times \mb{S}^4$, 
$\mb{S}^2\times\mb{S}^6$, and $\mb{S}^6\times \mb{S}^6$, respectively, and 
prove that none of these manifolds carry integrable almost complex structures.
\end{abstract}

\maketitle

\section{First: The Main Theorem}
As far back as 1951, using the Steenrod squares operations in algebraic 
topology, Borel and Serre proved that the only spheres that can
admit almost complex structures are $\mb{S}^2$, and $\mb{S}^6$ \cite{bs}.
For dimensional reasons, $\mb{S}^2$ admits an integrable structure, one 
induced by the existence of holomorphic coordinates on it. On the other hand,
using the purely imaginary Cayley numbers, and an embedding of $\mb{S}^6$ into
$\mb{S}^7\subset \mb{R}^8$, we can define an almost complex structure on
the 6-sphere, which is not integrable \cite{ehli} (the same method produces
almost complex structures on oriented hypersurfaces in 
$\mb{S}^7$ \cite{ca}, for instance, the products 
$\mb{S}^1(\sqrt{\frac{1}{6}})\times 
\mb{S}^5(\sqrt{\frac{5}{6}})$, $\mb{S}^2(\sqrt{\frac{1}{3}})\times 
\mb{S}^4(\sqrt{\frac{2}{3}})$, and
$\mb{S}^3(\sqrt{\frac{1}{2}})\times \mb{S}^3(\sqrt{\frac{1}{2}})$, 
of which the first and last, but not the second, are integrable). The 
question arose if the $6$-sphere, which carries this nonintegrable almost
complex structure, admits an integrable one.

The majority of the attempts to answer this question have failed to provide 
a definitive answer to it, as they do not elucidate how a homotopy 
between a plausible integrable structure, and the explicit one above on 
$\mb{S}^6$, varies as a function of the homotopy parameter, and have
been thus unable to decide if the plausible structure exists at all. We invoke
the Nash embedding theorem \cite{nash} to devise a way of overcoming
this obstacle: Since all Riemannian metrics on a closed manifold admit 
isometric embeddings into a sphere of large, but fixed dimension, we may  
identify all metrics on the manifold with their isometric embeddings into this 
fixed background, and by upgrading homotopies of almost complex structures to 
homotopies of almost Hermitian metrics associated to them in a natural way, 
we try to study the deformations of the former in terms of the extrinsic 
quantities associated to the isotopic deformations of the isometric embeddings 
of the latter. This reduces the issue in hand to the finding of a convenient 
path of such isometric embeddings of almost Hermitian metrics, which could 
then lead to its successful resolution.

The Yamabe problem on a Riemannian manifold $(M^n,g)$, and the $J$ Yamabe 
problem on an almost Hermitian manifold $(M,J,g)$, single out both the standard 
metric on $\mb{S}^6$, and the almost Hermitian pair that this metric and the 
octonionic almost complex structure above define together.
The first seeks a metric in the conformal class of $g$ of constant 
scalar curvature $s_g$, while the second seeks a metric in the conformal class 
of $g$ of constant $J$ scalar curvature $s_g^J$. They are solved by 
introducing the volume normalized total scalar, or $J$ scalar curvature 
functionals, whose infima are conformal invariants of the starting structures, 
bounded above by universal constants, and then proving that these 
conformal invariants are achieved by metrics in the class with the desired 
properties \cite{ya,au,tr,sc}, \cite{rss}. When the Yamabe conformal  
invariant is equal to the universal bound, and $M^n$ is the standard sphere, 
up to a conformal deformation, the solution to the Yamabe problem in the class
is the standard metric $g$ on $\mb{S}^n$, of trivial Weyl tensor $W_g$. In 
general, the Yamabe and $J$ Yamabe metrics on $(M,J,g)$ coincide if
we just have that $W_g(\omega_g^{\sharp}, \omega_g^{\sharp})=0$, 
$\omega_g$ the fundamental form, but unless $(M,g)$ is the standard 
sphere $(\mb{S}^6,g)$, these solutions do not achieve the universal bounds of 
their corresponding conformal invariants. In particular, the Yamabe and 
$J$ Yamabe metrics in a conformal class coincide when $(M,J,g)$ is locally 
conformally flat, and the solution to the $J$ Yamabe problem for the standard 
octonionic almost Hermitian structure $(\mb{S}^6,J,g)$ is, up to a conformal 
deformation, the standard metric $g$ also.

We can then view these abstractly defined metrics on $\mb{S}^{n}$, 
or $(\mb{S}^6, J)$, in terms of their Nash isometric embeddings into
the standard sphere $(\mb{S}^{\tn},\tg)$ of large dimension, and relate
their properties to special properties of the extrinsic 
quantities associated with their embeddings. We have that 
$s_g=n(n-1) + \| H_{f_g}\|^2 -\| \alpha_{f_g}\|^2$, where $H_{f_g}$ and 
$\alpha_{f_g}$ are the mean curvature vector and second fundamental form of the 
isometric embedding $f_g$, respectively, but if $g$ 
is a Yamabe metric in its conformal class $[g]$, the functions 
$\| H_{f_g}\|^2$ and $\| \alpha_{f_g}\|^2$ are both constants, and $f_g$ is 
a critical point of the total scalar curvature, and squared $L^2$-norms of 
$H_{f_g}$ and $\alpha_{f_g}$ functionals, under volume preserving conformal 
deformations of it, the last two of these properties characterizing metrics of 
constant scalar curvature in the class $[g]$. In particular, the 
standard metric on $\mb{S}^n$ has $\| \alpha_{f_g} \|^2 = 0$, and $s_g=n(n-1)$. 
Similarly, we have the general identity $s_g^J= 6 - 
\< \alpha_{f_g}(e_i,Je_j), \alpha_{f_g}(e_j, Je_i\>$ for the 
$J$ scalar curvature of a six dimensional almost Hermitian 
manifold $(M^6,J,g)$, but if $g$ is a $J$ Yamabe metric, the function 
$\< \alpha_{f_g}(e_i,Je_j), \alpha_{f_g}(e_j, Je_i\>$ is constant, and 
$f_g$ is a critical point of its squared $L^2$-norm under volume 
preserving conformal deformations of it, a property that characterizes the 
pairs $(J,g)$ of constant $s_g^J$. In particular, the 
standard $J$ Yamabe metric on $\mb{S}^6$ has $\< \alpha_{f_g}(e_i,Je_j), 
\alpha_{f_g}(e_j, Je_i)\>=0$, and $s_g^J=6$. 

We can now imagine that a sufficiently regular path of almost Hermitian pairs,
inducing a path of Yamabe metrics that begins at the conformal 
class of a plausible Hermitian Yamabe $(\mb{S}^6,J',g')$, and ends at the 
standard almost Hermitian Yamabe $(\mb{S}^6,J,g)$,  
contains enough information to determine if such a $J'$ exists, and can be 
deformed {\it continuously} to the latter Cayley structure. Since 
along the Palais' isotopic deformations $f_{g_t}$  
of the Yamabe metrics, we have an associated path 
$t\rightarrow \| H_{f_{g_t}}\|^2$ of constant functions that is as regular as 
the starting path of metrics, we can dilate these metrics conveniently along 
the path, and produce a new path of Yamabe metrics whose Nash embeddings are 
all minimal, absolute 
minimizers of $\Psi_{f_g}$, the squared $L^2$-norm of $H_{f_g}$ functional. 
We may then apply the gap theorem of Simons \cite{si,cdck,bla}, in 
its reinterpreted form as critical points of $\Psi_{f_g}$ \cite{rss2}, to 
complete the circle of ideas, and answer the question in the negative: We 
track the constant functions $\| \alpha_{f_{g}} \|^2$ along the path of 
isometric embeddings of the dilated Yamabe metrics with minimal embeddings, 
and show that they vary discontinuously as we move from the initial metric, 
whose conformal class is invariant under the integrable $J'$, to the latter 
one, whose conformal class is invariant under the octonionic $J$, so the 
former cannot exist.  

Using the definitions, and technical results described in \S2, whose proofs and 
references are provided in \S3, we execute this scheme, and prove now the 
following.

\begin{theorem}
The six dimensional sphere $\mb{S}^6$ does not carry integrable almost complex
structures.
\end{theorem}

{\it Proof}. Let $(\mb{S}^6,J,g)$ be the standard octonionic almost
Hermitian structure on $\mb{S}^6$. Suppose that $J'$ is an almost 
complex structure on $\mb{S}^6$ that defines the same orientation on it 
as that defined by the octonionic $J$. We show that $J'$ cannot be integrable.

\begin{enumerate}
\item[A)] Given a path $[0,1]\ni t\rightarrow J_t$ of almost
complex structures on $\mb{S}^6$ that begins at $J_0= J'$ and ends at  
$J_1 = J$, we consider the path of metrics $[0,1]\ni t \rightarrow g_t$ 
defined by
$$
g_t( \, \cdot \, , \, \cdot \, ) = \frac{1}{2}(
g( \, \cdot \, , \, \cdot \, )+g( J_t \, \cdot \, , J_t \, \cdot \, ))  \, ,
$$
and produce a smooth path $t\rightarrow  (J_t, g_t)$ of
almost Hermitian structures on $\mb{S}^6$, which begins at the almost Hermitian
structure $(J',g_0)$, and ends at the standard octonionic almost Hermitian
structure $(J,g)$. We obtain associated paths $t \rightarrow [g_t]$ and
$t\rightarrow (J_t,[g_t])$ of conformal classes, and almost Hermitian pairs
$t\rightarrow (J_t,[g_t])$, respectively.

\item[B)] On the conformal class $[g_t]$, there exists a metric $g_t^Y$ of
constant scalar curvature $s_{g_t^Y}$ that realizes the conformal invariant
$\lambda(\mb{S}^6,[g_t])$, the infimum of the Yamabe functional
$\lambda$ in (\ref{func}) over $[g_t]$, and a
metric $g^{JY}_t$ of constant $J_t$ scalar curvature $s_{g_t^{JY}}^{J_t}$
that realizes the conformal invariant $\lambda^{J_t}(\mb{S}^6,[g_t])$,
the infimum of the $J_t$ Yamabe functional $\lambda^J$ in (\ref{func}) over
$[g_t]$. These metrics can be normalized to have volume
$\omega_6$, the volume of the standard $(\mb{S}^6,g)$ itself.
Since $g$ is a solution of both, the Yamabe and $J$ Yamabe problem on
$(\mb{S}^6,J,g)$, by Lemma \ref{l44}, we can make choices of these metrics
$g_t^Y$ and $g_t^{JY}$ to produce (at least) $C^2$ paths
$t \rightarrow g_t^Y$, and $t \rightarrow g_t^{JY}$, respectively, that end
both at $g$ when $t=1$.

\item[C)] By the Nash isometric embedding theorem \cite{nash}, for a
sufficiently large codimension $p$, we may find paths of isometric embeddings
into the standard sphere background,
$$
f_{g^Y_t}: (\mb{S}^6,g_t^Y) \hookrightarrow  (\mb{S}^{\tn= 6 +p },\tilde{g})\, ,
$$
and
$$
f_{g^{JY}_t} : (\mb{S}^6,J_t,g_t^{JY}) \hookrightarrow
(\mb{S}^{\tn= 6 +p },\tilde{g})\, ,
$$
respectively. In terms of the extrinsic quantities of these embeddings,
by (\ref{sce}) and (\ref{sjip}), we have that
$$
s_{g_t^Y}=\frac{1}{\omega_6^{\frac{2}{6}}}\lambda(\mb{S}^6,[g_t])=
30-(\| \alpha_{f_{g^Y_t}} \|^2- \| H_{f_{g^{Y}_t}}\|^2)
$$
and
$$
s_{g_t^{JY}}^{J_t}=
\frac{1}{\omega_6^{\frac{2}{6}}} \lambda^{J_t}(\mb{S}^6,[g_t])
 =6-
\< \alpha^{J_t}_{f_{g_t^{JY}}} (e^t_i,J_te^t_j),\alpha^{J_t}_{f_{g_t^{JY}}}
(e^t_j,J_te^t_i)\> \, , 
$$
respectively. By Theorem \ref{thgr10}, $t\rightarrow \| H_{f_{g^{Y}_t}}\|^2$
is a path of constant functions.

\item[D)] Up to an isometry, the embeddings $f_{g_1^Y}$ and $f_{g_1^{JY}}$
coincide with each other, and equal the standard isometric embedding
$f_g: (\mb{S}^6,g) \hookrightarrow (\mb{S}^{6+p},\tilde{g})$
of the six sphere as a totally geodesic submanifold of
$(\mb{S}^{6+p},\tilde{g})$. For at $t=1$, $g_1^Y=g=g^{JY}_1$ is the
standard metric on $\mb{S}^6$, which achieves the conformal invariant
$\lambda(\mb{S}^6, [g])$, and which achieves the conformal
invariant $\lambda^{\tilde{J}}(\mb{S}^6, [g])$ for any almost complex structure
$\tilde{J}$ compatible with (the conformal class of) the standard metric $g$
also, in particular the octonionic $J$. We have that
$$
s_g=s_{g_1^Y}=\frac{1}{\omega_6^{\frac{2}{6}}} \lambda(\mb{S}^6,[g])=30 =
\frac{5}{\omega_6^{\frac{2}{6}}} \lambda^J(\mb{S}^6,[g])\, .
$$

\item[E)] Since the path $[g_t]\rightarrow g_t^Y$ is (at least) $C^2$,  
and the path $t \rightarrow s_{g_t^Y}$ of constant functions is 
(at least) continuous, 
by (\ref{sce}) the function $\| \alpha_{f_{g^{Y}_t}}\|^2-\| 
H_{f_{g^{Y}_t}}\|^2$ is a $t$ dependent constant, and the path
$t \rightarrow \| \alpha_{f_{g^{Y}_t}}\|^2-\| H_{f_{g^{Y}_t}}\|^2$ is 
continuous.  By Theorem \ref{th7}, the individual paths 
$t\rightarrow 
\| \alpha_{f_{g^{Y}_t}}\|^2$, and $t\rightarrow \| H_{f_{g^{Y}_t}}\|^2$ are 
$t$ dependent constants also, and for each $t$, the embedding $f_{g_t^Y}$ is 
a stationary point of the extrinsic functionals $\Psi_{f_g}(M)/\mu_g^{2/N}$, 
and $\Pi_{f_g}(M)/\mu_g^{2/N}$, $\Psi_{f_g}$ and $\Pi_{f_g}$ as in (\ref{tsc}),
in the space of volume preserving conformal deformations of $f_{g_t^Y}$, 
respectively.

\item[F)] Since $\mb{S}^6$ is oriented compatibly with the orientation of
$\mb{S}^{6+p}\subset \mb{R}^{n+p+1}$, we may write the mean curvature vector as 
$H_{f_{g^{Y}_t}}=h_t\nu_{f_{g_t^Y}}$, where $\nu_{f_{g_t^Y}}$ is a normal 
section of the normal bundle of $f_{g_t^Y}(\mb{S}^6)\hookrightarrow 
\mb{S}^{n+p}$, and $h_t$ is a $t$-dependent scalar.
If $T=T^{\tau}+T^{\nu}$ is the decomposition into tangential 
and normal components of the variational vector field of the family of 
isometric embeddings $f_{g^{Y}_t}$, we have that
$$
\frac{d\mu_{g_t^Y}}{dt}=\left( {\rm div}(T^{\tau}) - \< T^{\nu},
H_{f_{g^{Y}_t}}\>\right) d\mu_{g_t^Y}\, ,
$$
and since the path of metrics $t \rightarrow g_t^Y$ is volume preserving, it
follows that the component $T^{\nu}$ of $T$ is $L^2$-orthogonal to
$H_{f_{g^{Y}_t}}$. Since all normal directions are conformal, by (E) above, the 
functionals $\Psi_{f_{g^Y_t}}$ and $\Pi_{f_{g^Y_t}}$ are stationary in any 
normal direction $L^2_{g^Y_t}$ orthogonal to $\nu_{f_{g_t^Y}}$ also. On the 
other hand, by Theorem \ref{newth8}, the path of constant functions 
$t\rightarrow  h_t^2$ is as regular as the path $t\rightarrow g_t$ is, so the
conformally dilated metrics 
$\tilde{g}_t^Y=(1+h^2_t/n^2)g^Y_t$ are Yamabe metrics in $[g_t]$ whose 
corresponding family of isometric embeddings
$$
f_{\tilde{g}_t^Y}: (\mb{S}^6,(1+h^2_{t}/n^2)g_t^Y)
\rightarrow (\mb{S}^{6+p},\tilde{g})
$$ 
is minimal for all $t$, and the path of constant functions 
$t\rightarrow \| \alpha_{f_{\tilde{g}_t^Y}}\|^2$ is
continuous. These $f_{\tilde{g}_t^Y}$s are critical points of 
the functional
$$
f_g  \rightarrow  \Psi_{f_g}(M)
$$
along variations of $f_g$ in all normal directions, $\nu_{f_{g_t^Y}}$ 
unrestricted included \cite{gracie}. 
By the minimality of $f_{g^{Y}_1}$, we have that 
$\tilde{g}_1^Y=g^Y_1=g$, and $f_{\tilde{g}_1^Y}$ and $f_{g_1^Y}$ coincide,
up to isometries of the background sphere $(\mb{S}^{6+p},\tg)$. 

\item[G)] The functions $t \rightarrow \| \alpha_{f_{\tilde{g}_t^Y}}\|^2$, and
$t \rightarrow \| H_{f_{\tilde{g}^Y_t}}\|^2$, both vanish at $t=1$. By 
continuity, the lower and upper estimates in (\ref{es}) hold for a 
nontrivial open neighborhood of $1\in [0,1]$, and by Theorem \ref{mt}, we 
conclude that for $t$s in this neighborhood, 
$\| H_{f_{\tilde{g}^Y_t}}\|^2=0=\| \alpha_{f_{\tilde{g}^Y_t}}\|^2$,
$[g_{\tilde{g}^Y_t}]=[g_t^Y]=[g_t]=[g]$, and 
 modulo an isometry of the background sphere, 
$f_{\tilde{g}^Y_t}=f_{g_t^Y}$ is the standard 
embedding of $\mb{S}^6$ as a totally geodesic submanifold of $\mb{S}^{6+p}$.  

\item[H)] There are no orthogonal complex structures on the standard
sphere $(\mb{S}^6,g)$ \cite{bl}. Hence, if we assume that the path $J_t$
starts at an integrable $J'=J_0$, we can then conclude that
$[{\tilde{g}^Y_0}]=[g^Y_0]=[g_0]\neq [g]$. By Theorem \ref{t21},
$(\mb{S}^6,g_0)$ is not locally conformally flat, and the 
conformal invariants $\lambda(\mb{S}^6,
[g_0])$ and $\lambda^{J'}(\mb{S}^6,[g_0])$ achieved by $g_0^Y$ and $g_0^{JY}$,
respectively, are strictly less than their corresponding
universal bounds. Since the embeddings $f_{\tilde{g}^Y_t}$ are all minimal,
by Theorem \ref{mt}, for any $t$ for which $[g_t] \neq [g]$, the inequality  
$$
\| \alpha_{f_{\tilde{g}^Y_t}}\|^2 > n \geq \frac{np}{2p-1}  
$$
holds. Thus, the path $[0,1]\ni t
\rightarrow \| \alpha_{f_{\tilde{g}^Y_t}}\|^2$, which takes on a value 
greater than $6$ at $t=0$, and which is identically $0$ in a neighborhood of 
$1$, must have a discontinuity somewhere in between, a contradiction.
\end{enumerate}
\qed
\medskip

Quite more recently, by a topological argument that computes the Chern 
character of the product, and uses the integrality of the Chern 
character of a complex vector bundle and the Bott periodicity theorem, 
it has been proved that the products of even dimensional spheres carrying  
almost complex structures consists of just the cases 
$\mb{S}^2 \times \mb{S}^2$, $\mb{S}^2 \times \mb{S}^4$, $\mb{S}^2 \times 
\mb{S}^6$, and $\mb{S}^6 \times \mb{S}^6$, respectively \cite{dasu}. In \S4
we proceed as above to show that the last three 
of these manifolds do not carry integrable almost complex structures either,
our other main result. As we shall see then, the argument in the very last 
part of our proof there links prominently to the family of $J_t$ Yamabe 
metrics $g_t^{JY}$ in the proof, exposing the role that its alterego 
family here seems to play tacitly in the proof above.
\bigskip

\section{Second: Definitions and statements of technical results}
\label{second}

If $(M^n,g)$ is a Riemannian manifold, the Ricci and scalar curvature tensors
are defined by 
\begin{equation}
r_{g}(X,Y)={\rm trace}\, L\rightarrow R^g(L,X)Y \, , \label{rt}
\end{equation}
and 
\begin{equation}
s_g = {\rm trace}_g r_g \, , 
\end{equation}
respectively. Here, $R^g (X,Y)Z=(\nabla^g_{X}\nabla^g_{Y}- \nabla^g_{Y}
\nabla^g_{X}- \nabla^g_{[X,Y]})Z$ 
is the Riemann curvature tensor of $g$.  

Similarly, if $M^{n=2m}$ carries an almost complex structure $J$, and 
$(M,J,g)$ is an almost Hermitian manifold, the $J$-Ricci tensor, and $J$-scalar
curvature are defined by
\begin{equation}
r^{J}_g(X,Y)={\rm trace}\, L\rightarrow -J(R(L,X)JY) \, , \label{cr}
\end{equation}
and 
\begin{equation}
\label{Jsc}
s_g^J = {\rm trace}_g r_g^J \, , 
\end{equation}
respectively.
In contrast to $r_g$, $r^J_g$ is not always symmetric, and              
generally speaking, it is not $J$-invariant either. We have that 
$$
r_g^J(X,Y)=r_g^J(JY,JX) \, ,
$$
and if $(M,J)$ is a complex manifold of K\"{a}hler type, and the metric
$g$ is K\"{a}hler, then $r^J_g=r_g$.

We have the relation 
\begin{equation} \label{eqn3}
(n-1)s_g^J-s_g=2(n-1)W_g(\omega_g^{\sharp},\omega_g^{\sharp})\, , 
\end{equation}
$W_g$ the Weyl tensor, and $\omega_g^{\sharp}$ the fundamental form viewed
as a bivector \cite[Proposition 6]{rss}. By the conformal invariance of
$W_g$, the scalar tensor $s_g^J$ is independent of the choice of
$J$ compatible with $g$ used to define it.

If 
\begin{equation} \label{inem}
f_g: (M^n,g)\hookrightarrow (\tm,\tg)
\end{equation}
is an isometric embedding into a fixed Riemannian background $(\tm, \tg)$, and 
$\alpha:=\alpha_{f_g}$ and $H:=H_{f_g}$ are the second fundamental form and 
mean curvature vector of the embedding, the tensors $r_g$, and $s_g$, of $g$ 
relate to tensors of $\tg$, and the extrinsic quantities of the embedding by 
\begin{equation}
r_g(X,Y) =  
r_{\tg}(X,Y) - \sum_{i=n+1}^{\tn} \tg(R^{\tg}(\nu_i,X)Y,\nu_i) +\tg(H,\a(X,Y))-
 \sum_{i=1}^n \tg(\a(e_i,X),\a(e_i,Y)) \, ,
\end{equation}
and
\begin{equation}
s_g  = \sum K^{\tg}(e_i,e_j) +\tg(H,H)- \tg(\a,\a) \, ,
\label{sce}
\end{equation}
respectively. Here, $\{ e_i\}$ and $\{ \nu_i\}$ are orthonormal frames of 
the tangent space to the submanifold, and normal bundle, and $R^{\tg}$, 
$r_{\tg}$, and 
$K^{\tg}$ are the Riemann curvature tensor, Ricci tensor, and the sectional 
curvature of $\tg$, respectively. In what follows, we shall sometimes use the
summation convention, and write, for instance,
$\sum_{i,j}K^{\tg}(e_i,e_j)=\tg(R^{\tg}(e_i,e_j)e_j,e_i)$. We had used this 
convention already, in the preamble to Theorem 1.  

When in addition we have an almost complex structure $J$ on 
$M=M^{n=2m}$ and $(M,J,g)$ is an almost Hermitian manifold, there
are analogous relationships between $r^J_g$ and $s^J_g$ 
and the extrinsic quantities of the embedding given by   
\begin{equation} \label{ct}
r_g^J(Y,Z) = \tg(R^{\tg}(e_i,Y)JZ,Je_i)- \tg(\alpha(e_i,JZ),\alpha(Y,Je_i))
\, ,
\end{equation}
and 
\begin{equation} \label{sji}
s_g^J   =  \tg(R^{\tg}(e_i,e_j)Je_j,Je_i)-\tg(\alpha(e_i,Je_j),
\alpha(e_j,Je_i)) \, , 
\end{equation}
respectively. If we introduce the tensor 
\begin{equation} \label{aJ}
\a^J(X,Y)=\a(X,JY) \, , 
\end{equation}
whose symmetric and antisymmetric components are given by
\begin{equation} \label{den}
\a_{\pm}^{J}(X,Y)=\frac{1}{2}( \a(X,JY)\pm \a(Y,JX))\, , 
\end{equation}
we have that
\begin{equation} \label{sjip}
s_g^J   = 
\tg(R^{\tg}(e_i,e_j)Je_j,Je_i) + \| \a^J_{-}\|^2 - \| \a^J_{+}\|^2 \, ,
\end{equation}
where
$$
\begin{array}{rcl}
\| \a^J_{+}\|^2 & = & 
\frac{1}{2}( \< \alpha(e_i,Je_j),\alpha(e_i,Je_j)\> +
\< \alpha(e_i,Je_j),\alpha(e_j,Je_i))\>) \, , \vspace{1mm}  \\
\| \a^J_{-}\|^2 & = & 
\frac{1}{2} ( \< \alpha(e_i,Je_j),\alpha(e_i,Je_j)\> -
\< \alpha(e_i,Je_j),\alpha(e_j,Je_i))\>) \, . 
\end{array}
$$
The expressions (\ref{sce}) and (\ref{sji}) for $s_g$ and $s^J_g$ 
correspond to one another if we make $H$ 
and $\alpha$ correspond to $\a_{-}^J$ and $\a_{+}^J$, respectively. 

We set $N=2n/(n-2)$. The Yamabe and $J$ Yamabe functionals are  
\begin{equation} \label{func}
\lambda(g) = \frac{\displaystyle \int_{M}
s_{g} \,  d\mu _{g} }{
{\displaystyle \left( \int_{M}d\mu _{g}\right)^{2/N}}}\, ,
\quad \lambda^{J}(g)=\frac{\displaystyle \int_{M}
s^J_{g} \,  d\mu _{g} }{
{\displaystyle \left( \int_{M}d\mu _{g}\right) ^{2/N}}} \, ,
\end{equation}
respectively. The former is defined for any $(M^n,g)$, while the 
latter is defined when $(M^{n=2m},J,g)$ is an almost Hermitian structure.
The infima of these functionals over the conformal class 
$[g]$ of $g$ yield the conformal invariants 
$$
\begin{array}{rcl}
\lambda(M,[g]) & := & \inf_{g\in [g]} \lambda(g)  \hspace{1mm}\leq \hspace{1mm}
n(n-1) \omega_n^{\frac{2}{n}}\, , \\
\lambda^J(M,[g]) & := & \inf_{g\in [g]} \lambda^{J}(g)
\hspace{1mm}\leq \hspace{1mm}
n \omega_n^{\frac{2}{n}}\, .
\end{array}
$$
In the universal bounds of these invariants on the right above,  
$\omega_n$ is the volume of the standard unit sphere $\mb{S}^n\subset
\mb{R}^{n+1}$. 

\begin{lemma} \label{l44}
In the $C^2$-topology in the space of metrics, and its quotient topology 
in the space of conformal classes of metrics, the mappings 
$g \rightarrow \lambda(M,[g])$ and $[g]\rightarrow \lambda(M,[g])$ 
are continuous. With these topologies in the metric factors of the
product topology on the space of almost Hermitian pairs, the
mappings $(J,g) \rightarrow \lambda^J(M,[g])$ and $(J,[g])
\rightarrow \lambda^J(M,[g])$ are continuous. 
\end{lemma}

A Yamabe metric in $[g]$ is a metric $g^Y\in [g]$ that
realizes the conformal invariant $\lambda(M,[g])$. Its scalar curvature
is constant, and we have that
\begin{equation} \label{coy}
\lambda(M,[g])=s_{g^Y}\mu_{g^Y}(M)^{\frac{2}{n}} \, .
\end{equation}
Similarly, a $J$ Yamabe metric $g^{JY}$ on $(M,J,g)$ is one that
realizes the conformal invariant $\lambda^J(M,[g])$. Its 
$J$ scalar curvature is constant, and we then have that 
\begin{equation} \label{coay}
\lambda^J(M,[g])=s^J_{g^{JY}}\mu_{g^{JY}}(M)^{\frac{2}{n}} \, .
\end{equation}

Yamabe and $J$ Yamabe metrics always exist \cite{ya,au,tr,sc},
\cite{rss}. If $M\cong \mb{S}^n$, a Yamabe metric that achieves the 
universal bound is either a constant multiple of the standard metric
on $\mb{S}^n$, or its image under the action of a 
conformal diffeomorphism.
On the other hand, if $W_g(\omega_g^\sharp, \omega_g^\sharp)=0$, 
by (\ref{eqn3}), the conformal invariant that a Yamabe, and $J$ Yamabe metric
achieve are equal to each other (up to the factor $n-1$) 
\cite[Theorem 6.2]{rss}, and so, if $M\cong \mb{S}^n$ and the invariants
equal the universal upper bound, then $n=6$, and the $J$ Yamabe metric, 
which then can be taken to coincide 
with the Yamabe metric, is the metric on the standard six sphere, or an 
image of it under a conformal diffeomorphism, and 
the almost complex structure must be an almost complex structure compatible 
with it, for instance, the octonionic one introduced earlier. 

\begin{theorem} \label{t21}
Let $(M,J,g)$ be an almost Hermitian manifold. Then the solution to the
Yamabe and almost Hermitian Yamabe problem coincide if, and only if,
$$
W_g(\omega^\sharp_g, \omega^\sharp_g)=0 \, .
$$
\end{theorem}

By the Nash embedding theorem \cite{nash}, all Riemannian metrics on 
$M$ can be isometrically embedded into a standard sphere 
$(\mb{S}^{\tn},\tilde{g})$ of sufficiently large, but fixed dimension $\tn$. 
We thus identify a metric $g$ on $M$ with its isometric embedding $f_g$ into
this background sphere, and by (\ref{sce}), decompose the 
total scalar curvature functional as  \cite{gracie}
\begin{equation} \label{tsc}
\int_{f_g(M)} s_g d\mu_g = \Theta_{f_g}(M)+\Psi_{f_g}(M) - \Pi_{f_g}(M) :=
\Theta_{f_g}(M)-\mc{S}_{f_g}(M)\, , 
\end{equation}
where the expressions on the right are the extrinsic functionals  
$$
\begin{array}{ccl}
 \Pi_{f_g}(M) & = & 
{\displaystyle \int_M \| \alpha_{f_g}\|^2 d\mu_g } \, ,\vspace{1mm} \\
\Psi_{f_g}(M) & = & 
{\displaystyle \int_M \| H_{f_g}\|^2 d\mu_g } \, , \vspace{1mm} \\
 \Theta_{f_g}(M) & = 
& {\displaystyle \int_M \sum K^{\tg}(e_i,e_j) d\mu_g }\, . 
\end{array}
$$
We use this decomposition to analyze properties of Yamabe metrics in 
terms of properties of extrinsic quantities associated to their isometric
embeddings. 

If the scalar curvature $s_g$ of a metric $g$ on $M$ is a nonpositive 
constant, then $g$ is a Yamabe metric in its conformal class. In the
positive case, we have the following. 

\begin{theorem} \label{th4} 
If $g$ is a Riemmanian  metric of constant positive scalar curvature, and
there exists a minimal isometric embedding 
$f_g : (M,g) \hookrightarrow (\mb{S}^{\tn},\tilde{g})$ into the 
standard sphere, then $g$ is a Yamabe metric in its conformal class.  
\end{theorem}

\begin{corollary} \label{coam}
The standard product metric on  
$\mb{S}^{k}(\sqrt{k/n})\times \mb{S}^{n-k}(\sqrt{(n-k)/n})\subset \mb{S}^{n+1}$
is a Yamabe metric in its conformal class. The standard metrics on 
$\mb{P}^{n}(\mb{R}) \hookrightarrow \mb{P}^{n}(\mb{C})$ are Yamabe metrics in 
their respective conformal classes, and in their isometric minimal 
embeddings, the real projective space is the set of real points of 
complex projective space.  
\end{corollary}

{\it Proof}. All the manifolds 
$\mb{S}^{k}(\sqrt{k/n})\times \mb{S}^{n-k}(\sqrt{(n-k)/n})\subset \mb{S}^{n+1}$,
$1\leq k \leq n$, occur at the upper end of Simons' gap theorem 
for minimal embeddings into $\mb{S}^{n+p}$ satisfying $\| \alpha \|^2 \leq 
np/(2p-1)$ \cite[Theorem 5.3.2, Corollary 5.3.2]{si} 
\cite[Main Theorem]{cdck} \cite[Corollary 2]{bla}. Their linear embeddings  
into $\mb{S}^{n+1}$ are such that $\| \alpha \|^2=n$, so $s_g=n(n-2)$.
The projective space $\mb{P}^2(\mb{R}) \hookrightarrow \mb{S}^4$ also 
happens at this end of the gap theorem, minimally embedded into $(\mb{S}^4,\tg)$
with $\| \alpha \|^2 = 4/3$,  and $s_g=2/3$.   

Using eigenfunctions of the Laplacian associated to the first nonzero 
eigenvalue, we can define inductively isometric minimal embeddings of 
$\mb{P}^n(\mb{R})$, and $\mb{P}^n(\mb{C})$ (with their canonical metrics 
conveniently scaled) into the standard unit sphere, which present the former 
as the restriction of the latter to the set of real points \cite{sim}. 
(The novelty of this result lies in the method of proof, otherwise following
by the theorem that Einstein metrics are Yamabe metrics in their conformal 
classes \cite{oba,au}.)
\qed

Given an embedding $M^n \hookrightarrow \tm^{n+p}$ that is (at least) of class
$C^2$, we may write the mean
curvature vector as $H=h\nu_H$, where $\nu_H$ is a normal vector in the 
direction of $H$, and denote by $A_{\nu_H}$ and $\nabla^{\nu}$ the shape 
operator in the direction of $\nu_H$, and covariant derivative of the normal 
bundle, respectively. The local choice of direction $\nu_H$ may depend upon 
an orientation sign, as does the choice of $h$, but their product is globally 
well-defined, as are the quadratic expressions $\| \nabla^\nu \nu_H\|^2$ 
and $h^2=\| H\|^2$, respectively.

\begin{theorem} \label{thgr10}
If $n\geq 2$, the extrinsic function $\| H_{f_g}\|^2$ of any 
isometric embedding $f_g: (M^n,g)\hookrightarrow (\mb{S}^{n+p}, \tilde{g})$ is 
constant, and if $n \geq 3$, $f_g$ is a critical point of the functional 
$\Psi_{f_g}(M)/\mu_g^{2/N}$ in the space of volume preserving conformal 
deformations of the embedding.
\end{theorem}
\smallskip

\begin{theorem} \label{th7}
If $n\geq 3$, the metric $g$ of an isometric embedding 
$f_g : (M^n,g) \hookrightarrow (\mb{S}^{\tn},\tilde{g})$ is of
constant scalar curvature if, and only if,  the functions 
$\| \alpha_{f_g}\|^2$, and $\| H_{f_g}\|^2$ are constants, and
$f_g$ is a critical point of the extrinsic functionals 
$\Psi_{f_g}(M)/\mu_g^{2/N}$, and $\Pi_{f_g}(M)/\mu_g^{2/N}$ 
in the space of volume preserving conformal deformations of the embedding,  
respectively.
\end{theorem}
\medskip

\begin{theorem}\label{newth8}
Suppose that $M$ is orientable, and let $t\rightarrow g_t$ be a $C^k$ path
of metrics, $k\geq 2$, whose path $t\rightarrow [g_t]$ of classes has a lift   
to a path $t\rightarrow g_t^Y$ by Yamabe metrics $g_t^Y$,  
with isometric embeddings    
$$
[0,1]\ni t \rightarrow f_{g^Y_t} : (M^n,g_t^Y) \hookrightarrow
(\mb{S}^{\tn},\tilde{g}) \, ,
$$
that is minimal for at least one $t$. 
Then $t\rightarrow \| H_{f_{g^Y_t}}\|^2$ is a $C^k$ path of constant functions,
and the family $\tg_t^Y=(1+\| H_{f_{g^Y_t}}\|^2/n^2)g_t^{Y}$ of
 dilation deformations
of the $g_t^Y$s, is a path of Yamabe metrics whose associated family of
isometric embeddings
$$
[0,1]\ni t \rightarrow \tilde{f}_{\tg^Y_t} : 
(M^n,\left(1+\| H_{f_{g_t^Y}}\|^2/n^2
\right) g_t^Y) \hookrightarrow  (\mb{S}^{\tn},\tilde{g})
$$
is minimal for all $t$, and $\tilde{f}_{\tg^Y_t}=f_{g^{Y}_t}$ for any $t$ for 
which $\| H_{f_{g^Y_t}}\|=0$.   
\end{theorem}
\medskip

We consider embeddings, of the regularity indicated above, which satisfy the 
estimates
\begin{equation} \label{es}
\begin{array}{rcl}
{\displaystyle -\lambda \| H\|^2 -n} & \leq &
{\displaystyle {\rm trace}\, A_{\nu_H}^2-\| H\|^2 -\| \nabla^\nu \nu_H \|^2} \vspace{1mm} \\
 & \leq & {\displaystyle \| \a \|^2 -\| H\|^2 -\| \nabla^\nu \nu_H \|^2}
\leq {\displaystyle \frac{np}{2p-1}}
\end{array}
\end{equation}
for some constant $\lambda$. 
Notice that $\| A_{\nu_H}\|^2={\rm trace}\,
A_{\nu_H}^2$ is bounded above by
$\| \alpha \|^2$, and so the second of these inequalities is always true.
\smallskip

A normal stationary critical point of the functional $f_g \rightarrow
\Psi_{f_g}(M)$ is an isometric embedding 
$f_g: (M,g) \hookrightarrow (\tm,\tg)$ that is stationary in the
normal directions under deformations of the embedding.

\begin{theorem}\label{mt}
Suppose that $f_g: (M^n,g) \rightarrow (\mb{S}^{n+p},\tg)$, $n\geq 3$, 
is a normal
stationary critical point of the functional $f_g \rightarrow \Psi_{f_g}$. Then:
\begin{enumerate}
\item If the first inequality in {\rm (\ref{es})} 
holds for some constant $\lambda \in [0,1/2)$, then $f_g$ is a minimal 
embedding. In that case, if $0\leq \| \alpha_{f_g} \|^2 \leq np/(2p-1)$, we 
have that either 
\begin{enumerate}
\item $\| \alpha_{f_g}\|^2=0$, in which case $M$ lies in an equatorial sphere, 
or
\item $\| \alpha_{f_g}\|^2=np/(2p-1)$, in which case either $p=1$ and $M^n$ 
is any of 
the products $\mb{S}^{m}(\sqrt{m/n})\times \mb{S}^{n-m}(\sqrt{(n-m)/n})
\subset \mb{S}^{n+1}$, $1\leq m< n$, with $\| \alpha_{f_g} \|^2=n$,
or $n=p=2$ and $M^2$ is the real projective
plane embedded into $\mb{S}^4$ by the Veronese map, with $\| \alpha_{f_g}
\|^2=4/3$ and
scalar curvature $2/3$, all cases of metrics with nonnegative Ricci tensor, 
\end{enumerate}
while otherwise, the embedding is such that $\| \alpha_{f_g} \|^2>
n\geq np/(2p-1)$ 
at at least one point of $M$.
\item If the first inequality in {\rm (\ref{es})} fails to hold for any 
constant $\lambda \in [0,1/2)$, then $f_g$ cannot be minimal, and we have
that 
$$
0 \leq {\rm trace}\, A_{\nu_H}^2=
\frac{1}{2}\| H_{f_g}\|^2 + \| \nabla^{\nu}\nu_H\|^2 -n \leq 
\| \alpha_{f_g}\|^2 \, . 
$$ 
In particular, there is no $C^2$ deformation of normal stationary 
critical embeddings of this type into a minimal one.
\end{enumerate}
\end{theorem}
\medskip

\begin{example} \label{exgracie}
Given radii $r_1, r_2$, and an  $n \in \mb{N}$, we define the product 
Riemannian manifold
$$
M^{n}_k(r_1,r_2)= \mb{S}^k(r_1) \times \mb{S}^{n-k}(r_2) \subset
\mb{S}^{n+1}\left(\sqrt{r_1^2 + r_2^2}\right) \, .
$$
We set
$$
\mb{S}^{n,k} = M^n_k\left(\sqrt{\frac{k}{n}},\sqrt{\frac{n-k}{n}}\right)
\subset \mb{S}^{n+1} \, ,
$$
and, for convenience, name the case of even $n$ and odd $k$ separately,
setting
$$
C^{m,n} = M^{2m+2n+2}_{2m+1}\left(\sqrt{\frac{2m+1}{2m+2n+2}},\sqrt{
\frac{2n+1}{2m+2n+2}}\right)
\subset \mb{S}^{2m+2n+3} \, .
$$
The metrics on these manifolds are denoted by 
$g_{\mb{S}^{n,k}}$ and $g_{\mb{C}^{m,n}}$, respectively. 
As Riemannian manifolds,
$\mb{S}^{n,k}\simeq \mb{S}^{n,n-k}$, and with their natural complex structure
$J$, $(C^{m,n},J,g_{\mb{C}^{m,n}})$ is Hermitian \cite{caec}, and the Riemannian
isometry $C^{m,n}\rightarrow C^{n,m}$ reverses orientation.   

We have: 
\begin{enumerate} 
\item The linear isometric embedding $f_{g_{\mb{S}^{n,k}}}: (\mb{S}^{n,k},
g_{\mb{S}^{n,k}}) \hookrightarrow (\mb{S}^{n+1},\tg)$ is minimal and has
$\| \alpha_{f_{g_{\mb{S}^{n,k}}}}\|^2=n$. The scalar curvature of this
Yamabe metric is $n(n-2)$, and the Yamabe invariant of its conformal class is
$$
\lambda(\mb{S}^{n,k},[g_{\mb{S}^{n,k}}])=(n-2)(k^{\frac{k}{2}}(n-k)^{\frac{
n-k}{2}}\omega_k \omega_{n-k})^{\frac{2}{n}} \, .
$$
The sequence $k \rightarrow \lambda(\mb{S}^{n,k},[g_{\mb{S}^{n,k}}])$ is
decreasing in the range $1\leq k \leq [\frac{n}{2}]$. 

\item \label{two} The $J$-scalar curvature of the metric $g$ on 
$(C^{m,n},J,g_{\mb{C}^{m,n}})$ is 
$s^J_{g_{C^{m,n}}}=(2m/(2m+1)+ 2n/(2n+1))2(m+n+1)$, and we have that
$$
2(2m+2n+1)W_{g_{C^{m,n}}}(\omega_{g_{C^{m,n}}}^\sharp,
\omega_{g_{C^{m,n}}}^\sharp)= \frac{16mn(m+n+1)^2}{(2m+1)(2n+1)}\, .
$$   
so $W_{g_{C^{m,n}}}(\omega_{g_{C^{m,n}}}^\sharp,
\omega_{g_{C^{m,n}}}^\sharp)=0$ if either $m$, or $n$ is zero.
The almost complex Yamabe functional of this metric is
$$
\lambda^J(g_{C^{m,n}})=\left( \frac{2n}{2n+1}+\frac{2m}{2m+1}\right) \pi 
\left( \frac{4\pi}{n! m!} (2n+1)^{n+\frac{1}{2}}
(2m+1)^{m+\frac{1}{2}}\right)^{\frac{1}{n+m+1}} \, . 
$$
By Theorem \ref{t21}, this Yamabe metric is $J$-Yamabe also if
either $m$, or $n$ is zero. In any other case, the value of 
$\lambda^J(g_{C^{m,n}})$ above is strictly larger than the universal 
bound for the $J$-Yamabe problem, so this metric is not a $J$-Yamabe metric 
then.
\item \label{exgr3} Except in the cases of the standard $2$, and $6$ spheres, 
none of the Yamabe metrics on the almost Hermitian manifolds 
$\mb{S}^2\hookrightarrow 
\mb{S}^3$, $\mb{S}^2(\sqrt{1/2})\times \mb{S}^2(\sqrt{1/2})\hookrightarrow 
\mb{S}^5$, $\mb{S}^6\hookrightarrow \mb{S}^7$, 
$\mb{S}^2(\sqrt{1/3})\times \mb{S}^4(\sqrt{2/3})\hookrightarrow \mb{S}^7$, 
$\mb{S}^2(1/2)\times \mb{S}^6(\sqrt{3}/2)\hookrightarrow \mb{S}^9$, and  
$\mb{S}^6(\sqrt{1/2})\times \mb{S}^6(\sqrt{1/2})\hookrightarrow \mb{S}^{13}$ 
are $J$-Yamabe metrics also. For 
$\mb{S}^2(\sqrt{1/2})\times \mb{S}^2(\sqrt{1/2})$, we have
$s_g=8=s_g^J$, $\lambda([g])=\lambda^J(g)=8(2\pi)$, but the universal
bounds for the Yamabe and $J$-Yamabe problems are $12(
\frac{8\pi^2}{3})^{\frac{1}{2}}$, and $4( \frac{8\pi^2}{3})^{\frac{1}{2}}$, 
respectively.  For $\mb{S}^2(\sqrt{1/3})\times \mb{S}^4(\sqrt{2/3})$, we have 
that $s_g^J=12$, and $\lambda^J(g)=16\pi (\frac{2}{3})^{\frac{1}{3}}$, a value 
larger than the universal $J$-Yamabe bound $12\pi (\frac{2}{15})^{\frac{1}{3}}$.
For $\mb{S}^2(1/2)\times \mb{S}^6(\sqrt{3}/2)$ we have $s_g^J=16$, 
and $\lambda^J(g)=16 \pi (\frac{9}{20})^{\frac{1}{4}}$, a value larger than
the universal $J$-Yamabe bound $16\pi (\frac{2}{105})^{\frac{1}{4}}$.
For $\mb{S}^6(\sqrt{1/2})\times \mb{S}^6(\sqrt{1/2})$ we have
$s_g^J=24$, 
and $\lambda^J(g)=24 (\frac{2}{15})^{\frac{1}{3}}\pi$, a value larger than 
the universal $J$-Yamabe bound $8\pi(\frac{54}{385})^{\frac{1}{6}}$. 
\item For $r_n$ defined by $r_n^4=\left(\frac{n+1}{2}\right)^2(n-1)!$, we
view $\mb{P}^n(\mb{R})$, and $\mb{P}^n(\mb{C})$ with the metrics making the
$\mb{Z}/2$ and $\mb{S}^1$ fibrations 
$\mb{S}^n(r_n) \rightarrow \mb{P}^n(\mb{R})$ and $\mb{S}^{2n+1}(r_n) 
\rightarrow \mb{P}^n(\mb{C})$, respectively, Riemannian fibrations. 
With these scaled standard metrics on the projective spaces, we can then 
define minimal isometric embeddings $\iota^{\mb{R}}_n : 
\mb{P}^n(\mb{R}) \hookrightarrow \mb{S}^{\frac{n(n+3)}{2}-1}$ and
$\iota^{\mb{C}}_n: \mb{P}^n(\mb{C}) \hookrightarrow \mb{S}^{(n+1)^2-2}$ 
compatible with each other, in that the restriction of
$\iota^{\mb{C}}_n$ to the set of real points of $\mb{P}^n(\mb{C})$ equals
$\iota^{\mb{R}}_n$. The Yamabe invariants of these metrics are 
$$
\lambda(\mb{P}^n(\mb{R}),[g])= n(n-1) \left( \frac{\omega_n}{2}
\right)^\frac{2}{n}\, , \quad
\lambda(\mb{P}^n(\mb{C}),[g])= 4n(n+1) \left( \frac{\omega_{2n-1}}{2n}
\right)^\frac{1}{n}\, ,
$$
respectively. When $n\geq 2$, the real $2n$ dimensional manifolds 
$(\mb{P}^{2n}(\mb{R}),g)$, and $(\mb{P}^n(\mb{C}),g)$, which
embed minimally into $\mb{S}^{n(2n+3)-1}$ and
$\mb{S}^{(n+1)^2-2}$ with codimensions  
$n(2n+1)-1$ and $n^2-1$, respectively, have Yamabe 
invariants smaller than the Yamabe invariant 
of $(\mb{S}^{2n,n},g_{S^{n,k}})$, 
which embeds minimally into $\mb{S}^{2n+1}$ with codimension $1$. 
For $n=2$, the invariant of $(\mb{S}^{2n,n},[g_{S^{2n,n}}])$ is in between 
that of the projective spaces, with the one of the real projective space 
being always the smallest of the three.
\item For a sequence of values of $r\nearrow 1$, the product metric on 
the Calabi-Eckmann manifold $M_1^{2n}(r,\sqrt{1-r^2})=
\mb{S}^1(r) \times \mb{S}^{2n-1}(\sqrt{1-r^2}) \hookrightarrow \mb{S}^{2n+1}$ 
is Yamabe, and the corresponding sequence of Yamabe conformal invariants 
increases to the universal bound $2n(2n-1)\omega_{2n}^{\frac{1}{n}}$ 
\cite{sc2}. These manifolds are all locally conformally flat, and    
by Theorem \ref{t21}, their metrics are $J$ Yamabe metrics also. 
The manifold $(C^{0,n-1},J,g_{C^{0,n-1}})$ in (\ref{two}) above is the
first term of this sequence. This sequence of Hermitian 
manifolds becomes arbitrarily close conformally to the standard sphere 
$(\mb{S}^{2n},g)$, but the differing topologies, and complex structure of 
its terms, keep them apart from the limit.
\end{enumerate}
\end{example}

For the products of even dimensional spheres of Example \ref{exgracie} 
(\ref{exgr3}) above, with nonintegrable almost complex structures,
we have the following.

\begin{theorem} \label{nth9}
Let $(M^n,J,g)$ be any of the almost Hermitian Riemannian manifolds 
$(\mb{S}^{6,2},J,g_{\mb{S}^{6,2}})$, $(\mb{S}^{8,2},J,g_{\mb{S}^{8,2}})$, or 
$(\mb{S}^{12,6},J, g_{\mb{S}^{12,6}})$, and $f_g: (M,g) \hookrightarrow 
(\mb{S}^{n+1}, \tilde{g})$ be its standard minimal isometric embedding. If
$(-\varepsilon, \varepsilon) \ni t \rightarrow (\tilde{g}_t^Y,J_t)$ is a path 
of almost Hermitian deformations of $(g,J)$ by Yamabe metrics, with minimal 
associated family of isometric embeddings $f_{\tilde{g}_t^Y}: (M,\tilde{g}_t^Y) 
\hookrightarrow (\mb{S}^{n+p},\tilde{g})$, $f_{\tilde{g}_0^Y}=f_g$, 
and such that $[\tilde{g}_t^Y] \neq [g]$ for $t<0$, then $t\rightarrow 
f_{\tilde{g}_t^Y}$ is not $C^2$ continuous at $t=0$. In particular, the 
path $(-\varepsilon, \varepsilon)
 \ni t \rightarrow \| \alpha_{f_{\tilde{g}_t^Y}}\|^2$ changes discontinuously 
across $t=0$. 
\end{theorem}

\section{Third: Proofs of the technical results}

{\it Proof of Lemma \ref{l44}}. The mapping 
$g \rightarrow \lambda(M,[g])$ is continuous \cite[Proposition 7.2]{bebe}. 
By definition of the quotient topology, the mapping  
$[g] \rightarrow \lambda(M,[g])$ is continuous also.
The argument proving the first of these statements can be adapted readily 
to prove that $(J,g) \rightarrow \lambda^J(M,[g])$ is continuous. Then the
mapping $(J,[g]) \rightarrow \lambda^J(M,[g])$ is continuous also. \qed
\medskip

We observe next that if $(M,J,g)$ is an almost Hermitian manifold, then   
the skew Hermitian component of $r_g^J$ is a conformal invariant of 
$g$ \cite[Corollary 4.5]{rss}. This invariant encodes obstructions to the 
solutions of the Yamabe and $J$ Yamabe problems being the same.
When we pass to the trace of $r_g^J$, the said obstruction is 
captured by the function $W_g(\omega_g^\sharp, \omega_g^\sharp)$, and we then
see that the two Yamabe problems may have the same solution if
$(M,g)$ is locally conformally flat, but that this flatness is not a 
necessary condition for that to be the case.

Indeed, by (\ref{eqn3}), we have that
\begin{equation} \label{10}
\begin{array}{rcl}
(n-1)\lambda^J(M,[g])& \geq & \lambda(M,[g])+2(n-1)
\displaystyle{\inf_{g\in [g]}
\frac{1}{\mu_g(M)^{2/N}}
\int_M W_g(\omega_g^\sharp,\omega^\sharp_g)d\mu_g}
\, , \vspace{1mm} \\
\lambda(M,[g])& \geq & (n-1)\lambda^J(M,[g])+2(n-1)
\displaystyle{\inf_{g\in [g]} -\frac{1}{\mu_g(M)^{2/N}}
\int_M W_g(\omega_g^\sharp,\omega^\sharp_g)d\mu_g}
\, ,
\end{array}
\end{equation}
which allows us to extend \cite[Theorem 6.2]{rss} in the following way.

{\it Proof of Theorem \ref{t21}}.
If $W_g(\omega^\sharp_g,\omega^\sharp_g)=0$, by (\ref{eqn3})
we see that $(n-1)\lambda^J(M,[g])=\lambda(M,[g])$, and that the same metric
realizes the Yamabe and almost Hermitian Yamabe infimum invariant.

Conversely, if $g_0\in [g]$ is a Yamabe and $J$ Yamabe metric,
by (\ref{eqn3}), we conclude that
$W_{g_0}(\omega^\sharp_{g_0}, \omega^\sharp_{g_0})=c$, a constant. We write
the conformal invariants in terms of the conformal factor $\varphi\equiv 1$ 
relating $g_0$ to the Yamabe metric $\varphi^{\frac{4}{n-2}}g_0$. 
Since $\lambda^J(M,[g])=\lambda_N^{s_{g_0}^J}(1)$ and
$\lambda(M,[g])=\lambda_N^{s_{g_0}}(1)$, respectively, the first inequality
in (\ref{10}) implies that $c\leq 0$, while the second implies that $-c\leq 0$.
So $c=0$, and $W_g(\omega^\sharp_g,\omega^\sharp_g)=0$.
\qed
\medskip

We look back at $(M,g)$, and its isometric embedding $f_g$ (\ref{inem}) into the
background $(\tm,\tg)$. We develop suitable identities to study its isotopic 
conformal deformations.

Hence, let us suppose that
\begin{equation} \label{emc}
f_{g_t}: (M,g_t) \mapsto (\tm, \tg)
\end{equation}
is a path of isometric embeddings into the background associated to a 
path $t \rightarrow g_t$ of metrics deformations of $g=g_0$, with $f_{g_0}=f_g$.
For any $x\in M$, the trajectory of $f_g(x)=f_{g_0}(x) \in \tm$ is
given by the path $t\rightarrow f_{g_t}(x)$. Under mild hypothesis on
$(\tm,\tg)$, which hold for any spaceform background, we can apply the 
Palais' isotopic extension theorem \cite{pal}, and
obtain a smooth one parameter family of diffeomorphisms
$$
F_t: \tm \rightarrow \tm 
$$
such that
$$
F_t(f_g(x))=f_{g_t}(x)\, . 
$$
The pullback tensor $F_t^* \tg$ is just the background metric on $\tm$
acted on by the diffeomorphism $F_t$s. Thus, since all the metrics
are realized by the metric on the background manifold, we can 
use pull-back and restrictions to obtain a path of diffeomorphisms  
$$
F_t\mid_{f_g(M)}: (f_g(M),\tg) \rightarrow (f_{g_t}(M), \tg)\hookrightarrow
(\tm, \tg)\, ,
$$
identifying the path of isometric embeddings deformations $f_{g_t}$ with 
a path of equally regular isotopic deformations of the submanifold $f_g(M)$ in 
$\tm$.

Suppose now that the path (\ref{emc}) corresponds to a path of  
conformally related metrics
\begin{equation} \label{cd}
[0,1]\ni t \rightarrow g_t =e^{2\psi (t)}g\, .  
\end{equation}
Thus, $f_{g_t}$ is the associated family of isometric embeddings 
conformally deforming the initial embedding $f_{g}$. 
The pullback tensor $F_t^* \tg$ is just the metric $\tg$ on $\tm$, 
and by construction, we have that
$$
F_t^* (\tg\mid_{f_{g_t}(M)}) = e^{2 u(t)(f_g(\, \cdot \,))} \tg \mid_{f_g(M)} =
e^{2 u(t)(f_g(\, \cdot \,))} \tg\mid_{f_g(M)}  \, ,
$$
where the conformal factor $e^{2 u(t)}$ and that in (\ref{cd}) are related 
to each other by $e^{2\psi(t)(\, \cdot \, )}=
e^{2u(t)\circ f_g(\, \cdot \, )}$.  
We extend $u(t)$ conveniently to a function defined on the whole of
$\tm$, and view the family (\ref{emc}) as the family of conformal isometric 
embeddings 
\begin{equation}\label{embt}
f_{g_t}: (M,e^{2u(t)\circ f_g}g) \hookrightarrow 
(f_{g}(M),e^{2u(t)\circ f_g}g) \hookrightarrow (\tm ,\tg)
\end{equation}
of the fixed submanifold $f_g(M)\hookrightarrow \tm$ with a varying 
conformal metric on it, deforming $f_g(M)$ by conformal isotopies in $\tm$.   
We compute the intrinsic and extrinsic quantities in (\ref{sce})
associated to $f_{g_t}$, and relate them all to those of $f_g$, the 
initial isometric embedding. 
For notational convenience here, we set $\tilde{g}_t=e^{2u(t)}\tilde{g}$.

The curvature tensors of $\tilde{g}_t$ and $\tilde{g}$ are related to each
other by
$$ 
R^{\tilde{g}_t}  = e^{2u}(R^{\tilde{g}}+\tg \kuno
 (\nabla^{\tg} du-du\circ du + \frac{1}{2}|du|^2 \tg ))\, , 
$$
where $p \bigcirc \hspace{-3.40mm} \wedge \hspace{1.4mm} q$ 
is the Kulkarni-Nomizu product of the symmetric 2-tensors $p$ and $q$.
If $e_1, \ldots, e_n$ is a $\tilde{g}$-orthonormal tangent frame on $f_g(M)$, 
then $e_1^t=e^{-u(t)}e_1, \ldots, e_n^t=e^{-u(t)}e_n$ is an orthonormal 
tangent frame for $f_{g_t}(M)$, and we have that
$$ 
\sum_{i,j} K^{\tilde{g}_t}(e_i^t, e_j^t)= e^{-2u(t)}(
\sum_{i,j} K^{\tilde{g}}(e_i,e_j) -2(n-1) 
 {\rm trace}_{f_g(M)} p) \, ,
$$ 
where $p=\nabla^{\tg} du-du\circ du + \frac{1}{2}|du|^2 \tg$. If $\tau$ and
$\nu$ denote the tangential and normal components, respectively, we have that
\begin{equation} \label{eq2}
{\rm trace}_{f_g(M)} p=
{\rm div}_{f_g(M)} (\nabla^{\tg} u)^\tau- \tg(H_{f_g}, (\nabla^{\tg}
u)^\nu)-|du^\tau|_{\tg}^2+\frac{n}{2}|du|_{\tg}^2 \, , 
\end{equation}
and so
\begin{equation} \label{eq1}
\sum_{i,j} K^{\tilde{g}_t}(e_i^t, e_j^t)= e^{-2u(t)}(
\sum_{i,j} K^{\tilde{g}}(e_i,e_j) -2(n-1)(
{\rm div}_{f(M)} (\nabla^{\tg} u)^\tau- \tg(H_{f_g}, (\nabla^{\tg}
u)^\nu)-|du^\tau|_{\tg}^2+\frac{n}{2}|du|_{\tg}^2 )) \, .
\end{equation}

On the other hand, since the second fundamental forms of $f_{g_t}$ and 
$f_g$ are related to each other by 
$$
\alpha_{f_{g_t}}(X,Y)  = \alpha_{f_g}(X,Y)-\tg(X,Y)(\nabla u)^\nu\, ,
$$
and we obtain that 
\begin{eqnarray}
\| H_{f_{g_t}} \|^2 & = & e^{-2u(t)}(\| H_{f_g}\|^2 
- 2n\tg(H_{f_g},\nabla^{\tg}u ^\nu)
+n^2 \tg(\nabla^{\tg}u ^\nu, \nabla^{\tg} u^\nu)) \, , \label{eq3}\\
\| \alpha_{f_{g_t}}\|^2 & = & e^{-2u(t)}(\| \alpha_{f_g}\|^2 - 2\tg(H_{f_g},
\nabla^{\tg}u^\nu) +n \tg(\nabla^{\tg}u ^\nu, \nabla^{\tg} u^\nu)) \, .
\label{eq4} 
\end{eqnarray}

Relations (\ref{eq1}), (\ref{eq3}), and (\ref{eq4}) yield
\begin{equation} \label{eq5}
s_{\tg_t}= e^{-2u(t)}\left( s_{g}-2(n-1){\rm div}_{f_g(M)}(\nabla^{\tg}u)^\tau
-(n-1)(n-2)\tg(\nabla^{\tg}u^{\tau},\nabla^{\tg}u^\tau)  \right) \, , 
\end{equation}
as is to be expected.

Similarly, if $J$ is a Hermitian almost complex structure for the conformal
class of $g$, 
\begin{equation} \label{eq6}
s^J_{\tg_t}= e^{-2u(t)}\left( s^J_{g}-2{\rm div}_{f_g(M)}(\nabla^{\tg}u)^\tau
-(n-2)\tg(\nabla^{\tg}u^{\tau},\nabla^{\tg}u^\tau)  \right) \, .
\end{equation}

By the Nash isometric embedding theorem \cite{nash}, we may take the standard
sphere $(\mb{S}^{\tn},\tg)$ of sufficiently large dimension $\tn=\tn(n)$ to
play the role of the background Riemannian manifold in the construction 
above. There is an optimal choice of $\tn$ that depends on the dimension 
$n$ of $M$ only. By identifying metrics, and their deformations, with their 
isometric embeddings into this background, and their corresponding Palais 
isotopic deformations as above, we may proceed to characterize metrics of 
constant 
scalar curvature on $M$ by certain properties of the extrinsic quantities of 
their associated embeddings.

In terms of the extrinsic functionals defined in (\ref{tsc}), the
Yamabe functional (\ref{func}) of the conformally related metrics 
$\tilde{g}_t=e^{2u(t)}\tilde{g}$, with associated family of  
isometric embeddings (\ref{emc}), decomposes as the linear combination 
\begin{equation} \label{lu1}
\lambda(\tilde{g}_t) = \frac{1}{\mu^{\frac{2}{N}}_{\tilde{g}_t}}
\Theta_{f_{{g}_t}}(M) -
\frac{1}{\mu^{\frac{2}{N}}_{\tilde{g}_t}} \mc{S}_{f_{{g}_t}}(M) \, .
\end{equation}
In calculating the derivatives of these functionals, and finding the
weak equation that their critical points satisfy, the
set of all paths $t\rightarrow u(t)$ defining the conformal deformation in the
class plays the role of the set of test 
functions.
 
If the background is the standard sphere, whose sectional curvature is
$1$, the exterior scalar curvature $\sum_{i,j} K^{\tilde{g}_t}(e_i^t,e_j^t)$ 
is the constant $n(n-1)$, and by (\ref{eq1}) we conclude that 
\begin{equation} \label{lu2}
n(n-1)=e^{-2u(t)}(n(n-1)+(n-1)(2\Delta^{\tilde{g}}u -(n-2)
\tilde{g}( \nabla^{\tilde{g}} u^\tau, \nabla^{\tilde{g}} u^\tau)
+\tilde{g}(2H_{f_g} -n\nabla^{\tilde{g}} u^\nu, \nabla^{\tg} u^\nu )
)) \, .
\end{equation}
On the other hand, if the conformal deformations of the metrics are volume 
preserving, the constant exterior scalar curvature implies that
$\Theta_{f_{\tg_t}}/\mu^{\frac{2}{N}}_{\tg_t}$ is constant as well, and by
(\ref{eq1}), we conclude that  
\begin{equation}\label{basic}
0=\frac{d}{dt} \left( \frac{1}{\mu^{\frac{2}{N}}_{\tilde{g}_t}}\int 
\sum K^{\tilde{g}_t}( e_i^t , e_j^t ) d\mu_{\tilde{g}_t} \right) =
\frac{-2(n-1)}{\mu^{\frac{2}{N}}_{\tilde{g}_t}}\int_{f_g(M)} 
e^{(n-2)u(t)}\tilde{g}(n\nabla^{\tilde{g}}u^\nu -H_{f_g},
\nabla^{\tilde{g}} \dot{u}^\nu ) d\mu_{\tilde{g}} 
\end{equation} 
holds for all $t$. Thus, along arbitrary paths of volume preserving conformal 
deformations of the initial embedding $f_g$, the curve $u(t)$, 
which always satisfies the strong identity (\ref{lu2}), must have a velocity
$\dot{u}$ that satisfies this weak equation (\ref{basic}) as well. (Suitably 
accounting for the different curvatures, these statements remain equally true 
for isometric embeddings of volume preserving conformally related metrics 
into any space form background.) We exploit this result time and time again 
to compute the variational equation of functionals of worth in our problem.  

By (\ref{eq3}) and (\ref{eq4}), (\ref{basic}) implies that 
\begin{equation} \label{gr31}
\begin{array}{rcl}
{\displaystyle \frac{d}{dt} \left(\frac{1}{\mu^{\frac{2}{N}}_{\tilde{g}_t}}
\mc{S}_{\tilde{g}_t} \right)} \hspace{-3mm} & = & \hspace{-3mm}
 {\displaystyle \frac{n-2}{
\mu^{\frac{2}{N}}_{\tilde{g}_t}}\int \hspace{-1mm} \dot{u}\left(
\| \alpha_{f_{{g}_t}}\|^2 \!\!  -\| H_{f_{{g}_t}}\|^2 \hspace{-0.9mm} - 
\hspace{-0.9mm}
\frac{\mc{S}_{\tilde{g}_t}}{ \mu_{\tilde{g}_t}}\right)
d\mu_{\tilde{g}_t} - \frac{2(n-1)}{ \mu^{\frac{2}{N}}_{\tilde{g}_t}}
\int e^{(n-2)u(t)}\tilde{g}(n\nabla^{\tilde{g}}u^\nu \! -H_{f_g},
\nabla^{\tg} \dot{u}^\nu ) d\mu_{\tilde{g}} } \vspace{1mm} \\
& = & {\displaystyle  \frac{n-2}{
\mu^{\frac{2}{N}}_{\tilde{g}_t}}\int \dot{u}\left(
\| \alpha_{f_{{g}_t}}\|^2-\| H_{f_{{g}_t}}\|^2-
\frac{\mc{S}_{\tilde{g}_t}}{ \mu_{\tilde{g}_t}}\right)
d\mu_{\tilde{g}_t} } \, , 
\end{array}
\end{equation}
so embeddings for which $\| \alpha_{f_g}\|^2-\| H_{f_g}\|^2$ is constant are
precisely the critical points of 
$\mc{S}_{\tilde{g}_t}/\mu^{\frac{2}{N}}_{\tilde{g}_t}$
among volume preserving conformal deformations, and when that is so, the
critical point equation is equivalent to the vanishing of the expression 
in the right side of (\ref{basic}) at the corresponding $t$, and for all 
possible $\dot{u}$s.    
\medskip

{\it Proof of Theorem \ref{th4}}.  We let $g'$ be a Yamabe metric in 
the conformal class of $g$ of the same volume as that of $g$, and consider 
arbitrary paths $g_t$ of volume preserving conformal deformations (\ref{cd}) 
that start at $g$ when $t=0$, and end at $g'=g_1$ when $t=1$. We apply the
construction above, and produce a path of 
volume preserving conformal deformations (\ref{embt}) of the embedding 
$f_g$ to an isometric embedding $f_{g'}$ of $g'$. We let $u(t)$ be the   
associated path of scalar functions such that 
$\tilde{g}_t =e^{2u(t)}\tilde{g}$.

Since the scalar curvature of $\tilde{g}_t$ is
constant both when $t=0$, and $t=1$, by (\ref{sce}), and the variational 
expression (\ref{gr31}), it follows that
$\left( \frac{1}{\mu^{\frac{2}{N}}_{\tilde{g}_t}} \mc{S}_{\tilde{g}_t} \right)$ 
has a critical point at these two values of $t$.
Since the embedding $f_g$ is minimal, by (\ref{eq3}) and (\ref{eq4}) we 
conclude that the constant $\| \alpha_{f_{{g'}}}\|^2-\| H_{f_{{g'}}}\|^2=
\| \alpha_{f_{{g}_1}}\|^2-\| H_{f_{{g}_1}}\|^2$ is equal to
\begin{equation} \label{lu3}
\| \alpha_{f_{{g'}}}\|^2-\| H_{f_{{g'}}}\|^2 = e^{-2u(1)}( 
\| \alpha_{f_g}\|^2 +n(1-n)\tilde{g}(\nabla^{\tilde{g}}u^\nu, 
\nabla^{\tilde{g}}u^\nu)\mid_{t=1} ) \, . 
\end{equation}

The quotient of this, and  (\ref{lu2}) evaluated at $t=1$, after some 
simplifications, yields that 
$$
n s_{g'} \tilde{g}(\nabla^{\tilde{g}}u^\nu,   
\nabla^{\tilde{g}}u^\nu)\mid_{t=1} =n(s_{g'}- s_{g}) 
+(\| H_{f_{{g'}}}\|^2 -\| \alpha_{f_{{g'}}}\|^2)
(2\Delta^{\tilde{g}}u -(n-2)\tilde{g}(\nabla^{\tilde{g}}u^\tau, 
\nabla^{\tilde{g}}u^\tau)\mid_{t=1} \, , 
$$
which by integration with respect to the measure 
$e^{\frac{(n-2)u(1)}{2}}d\mu_{\tilde{g}}$ produces
$$
0\leq n s_{g'} \int e^{\frac{(n-2)u(1)}{2}} \tilde{g}(\nabla^{\tilde{g}}u^\nu, 
\nabla^{\tilde{g}}u^\nu)\mid_{t=1} d\mu_{\tilde{g}} =
n(s_{g'}- s_{g}) \int e^{\frac{(n-2)u(1)}{2}}  d\mu_{\tilde{g}} \leq 0 \, .
$$ 
Thus, $\nabla^{\tilde{g}}u^{\nu}\mid_{t=1}=0$, and $s_{g'}=s_{g}$. By 
(\ref{lu3}), it follows that the function $u(1)$ is constant. Since $g'$ and 
$g$ have the same volume, this constant is zero. So $g'=g$. \qed  
\medskip

The identities (\ref{basic}), (\ref{eq3}) and (\ref{eq4}) imply that 
\begin{equation} \label{gr32}
\begin{array}{rcl}
{\displaystyle \frac{d}{dt} \left(\frac{1}{\mu^{\frac{2}{N}}_{\tilde{g}_t}}
{\Psi}_{\tilde{g}_t} \right)} \! \! & = & \! \!  
{\displaystyle  \frac{n-2}{
\mu^{\frac{2}{N}}_{\tilde{g}_t}}\int \dot{u}\left(
\| H_{f_{{g}_t}}\|^2-
\frac{{\Psi}_{\tilde{g}_t}}{ \mu_{\tilde{g}_t}}\right)
d\mu_{\tilde{g}_t} 
+
\frac{2n}{ \mu^{\frac{2}{N}}_{\tilde{g}_t}}
\int e^{(n-2)u(t)}\tilde{g}(n\nabla^{\tilde{g}}u^\nu \! -H_{f_g},
\nabla^{\tg} \dot{u}^\nu ) d\mu_{\tilde{g}}
} \\ & = & \! \!
{\displaystyle  \frac{n-2}{
\mu^{\frac{2}{N}}_{\tilde{g}_t}}\int \dot{u}\left(
\| H_{f_{{g}_t}}\|^2-
\frac{{\Psi}_{\tilde{g}_t}}{ \mu_{\tilde{g}_t}}\right)
d\mu_{\tilde{g}_t}
} \, , \\
{\displaystyle \frac{d}{dt} \left(\frac{1}{\mu^{\frac{2}{N}}_{\tilde{g}_t}}
{\Pi}_{\tilde{g}_t} \right)} \! \! & = & \! \!   
{\displaystyle  \frac{n-2}{
\mu^{\frac{2}{N}}_{\tilde{g}_t}}\int \dot{u}\left(
\| \alpha_{f_{{g}_t}}\|^2-
\frac{{\Pi}_{\tilde{g}_t}}{ \mu_{\tilde{g}_t}}\right)
d\mu_{\tilde{g}_t} 
+ 
\frac{2}{ \mu^{\frac{2}{N}}_{\tilde{g}_t}}
\int e^{(n-2)u(t)}\tilde{g}(n\nabla^{\tilde{g}}u^\nu \! -H_{f_g},
\nabla^{\tg} \dot{u}^\nu ) d\mu_{\tilde{g}}
} \\  
& = & \! \!
{\displaystyle  \frac{n-2}{
\mu^{\frac{2}{N}}_{\tilde{g}_t}}\int \dot{u}\left(
\| \alpha_{f_{{g}_t}}\|^2-
\frac{{\Pi}_{\tilde{g}_t}}{ \mu_{\tilde{g}_t}}\right)
d\mu_{\tilde{g}_t} 
} \, ,
\end{array}
\end{equation}
so embeddings for which $\| H_{f_g}\|^2$ and 
$\| \alpha_{f_g}\|^2$ are constant are precisely the critical points of
the extrinsic functionals 
$\Psi_{\tilde{g}_t}/\mu^{\frac{2}{N}}_{\tilde{g}_t}$, and
$\Pi_{\tilde{g}_t}/\mu^{\frac{2}{N}}_{\tilde{g}_t}$, respectively,  
in the space of volume preserving conformal deformations of the embedding. 
The critical points of the three extrinsic functionals in the combination 
(\ref{tsc}), which normalized by the volume have linear combination   
 (\ref{lu1}) equal to the Yamabe functional, are characterized as the 
critical points of the Yamabe functional itself: Their densities are 
constant functions, and when this is the case, their critical point equations
are equivalent to the vanishing of the right side of (\ref{basic}) at the
corresponding $t$, for all possible $\dot{u}$s. But the particular 
curvature properties of the sphere background brings about an 
additional affinity between the first two of these extrinsic functionals, 
which we unfold next. 
\smallskip

Given an isometric embedding $f_g: (M,g)\hookrightarrow (\mb{S}^{n+p}, 
\tilde{g})$, let us consider a local orthonormal frame $\{e_i\}_{i=1}^n$ of 
$TM$ valid in some neighborhood of $x\in f_g(M)$.  By applying Euclidean 
rotations that take a vector $v$ into a vector $w$ while leaving the orthogonal 
complement that these two vectors span unchanged, we find an Euclidean rotation
that transforms the standard ordered basis $\{c_1, \ldots, c_{n+p+1}\}$ for 
$T_x\mb{R}^{n+p+1}$ into an orthonormal ordered basis of the form 
$\{ e_1, \ldots, e_n, d_{n+1}, \ldots, d_{n+p+1}\}$, where the $d_j$s 
constitute an ordered orthonormal basis of the normal bundle of $f_g(M)$ at $x$
 in the ambient Euclidean space. If necessary, we may apply an additional 
Euclidean rotation to these normal vectors that transforms $d_{n+p+1}$ into 
the vector $x$ itself, thought as a unit tangent radial vector in 
$\mb{R}^{n+p+1}$, while leaving the orthogonal complement of their span fixed. 
We rename the frame so obtained as $\{ e_1, \ldots, e_n, d_{n+1}, \ldots, 
d_{n+p}, x\}$. 
The mean curvature vector $H_{f_g}$ at the point $x$ lies in the 
span of $\{ d_{n+1}, \ldots, d_{n+p}\}$, vectors that are a  
basis of $\nu_x(f_g(M))$.  If this vector is nonzero, we apply one more 
Euclidean rotation of the type being used that transforms $d_{n+1}$ into 
a unit vector $\nu_H$ in its direction; otherwise, we 
set $\nu_H=d_{n+1}$. The constructed resulting basis $\{ d_{n+1}=\nu_H, \ldots,
 d_{n+p},x\}$ of the normal bundle depends on the  choice of 
$\{ e_1, \ldots, e_n\}$ only by an orientation sign of this
frame, and varies smoothly with $x$ in the neighborhood where the frame is 
valid. Hence, if $M$ is orientable, by patching together local constructions, 
we obtain a globally defined normal vector field $\nu_H$, and a suitable scalar
function $h$, such that $H_{f_g}=h\nu_H$, with a 
line splitting off in the normal bundle of $f_g(M) \hookrightarrow 
\mb{S}^{n+p}$. If $M$ is not orientable, the vector 
field $\nu_H$ is only defined locally up to the choice of a direction, and
the scalar function $h$ is correspondingly defined up to the choice of a sign,
but the product $H=h \nu_H$ is well defined globally, although the lines   
defined by $\nu_H$ at each $x$ do not split off globally as a line bundle 
summand of the normal bundle of $f_g(M)$ in $\mb{S}^{n+p}$. 

Notice that if the background is a Riemannian manifold $(\tm,\tg)$ other than
the standard sphere, we can likewise produce a decomposition of the form
$H_{f_g}=h_{f_g}\nu_{f_g}$ for the mean curvature vector of an isometric 
embedding $f_g : (M,g) \rightarrow (\tm, \tg)$, with the factors of properties
similar to the ones they have in the case above. Indeed, we may use the Nash 
embedding theorem to get an isometric embedding $f_{(\tm,\tg)}: (\tm,\tg) 
\rightarrow (\mb{S}^{\tn}, \tilde{g})$ of this new background into a 
standard sphere of sufficiently large dimension, and apply the argument above 
iteratively to the composition 
$f_{\tm,\tg} \circ f_g$ to draw the desired conclusion in this more general 
context.
\smallskip

{\it Proof of Theorem \ref{thgr10}.} Suppose first that $M$ is oriented
compatibly with the orientation of the ambient space sphere $\mb{S}^{n+p}
\subset \mb{R}^{n+p+1}$. We consider the normal vector field $\nu_H$ above
such that $H_{f_g}=h \nu_H$, $h$ a real valued function on $f_g(M)$. 

Let $\pi_g$ be the $L^2$ projection operator onto the constant functions. 
At a point on $f_g(M)$, we define $t$ to be the arclength parameter of the 
geodesic in $(\mb{S}^{n+p},\tg)$ emanating from the submanifold in the 
direction of $\nu_H$ at the point. This procedure defines $t$ as a scalar 
function on $f_g(M)$ ranging in $(-a,a)$ for some sufficiently small positive
$a$.  If $\mu_{t}$ is the volume of the conformally related metric 
$e^{2(h-\pi_g h)t}g$, we rescale it to define the path of volume preserving
conformally related metrics $g_t=e^{2\left( (h-\pi_g h)t + 
\ln{\mu_g}-\ln{\mu_{t}}\right)}=e^{2 \psi(t)}$. 

Since the sphere $(\mb{S}^{n+p},\tg)$ is of sufficiently large dimension, and
so it carries isometric embeddings of all metrics on $M$, by the Nash isometric 
embedding theorem, we produce a family of volume preserving isometric 
embeddings $f_{g_t}: (M,g_t)\rightarrow (\mb{S}^{n+p},\tg)$ that are conformal 
deformations of $f_g$. If we view the metric on $f_{g_t}(M)$ as 
$e^{2u(t)}\tilde{g}$ in the conformal deformation (\ref{embt}) of $f_{g}$, by 
construction, $u(0) = 0$, ${\displaystyle \frac{d}{dt}}u 
\mid_{t=0}=(h-\pi_g h)$, and $\nabla^{\tilde{g}} \dot{u}^{\nu}\mid_{t=0}=
(h-\pi_g h)\nu_H$, respectively.
 
Since the sphere $(\mb{S}^{n+p},\tg)$ has constant sectional curvature, and
$n\geq 2$, the volume preserving condition is equivalent to the vanishing of 
the expression on the right side of (\ref{basic}) for all $t$s for which
$u(t)$ is defined, in particular, $t=0$. Thus,   
$$
\begin{array}{rcl} 
0 & = &  {\displaystyle \int_{f_g(M)}
e^{(n-2)u(t)}\tilde{g}(n\nabla^{\tilde{g}}u^\nu -H_{f_g},
\nabla^{\tilde{g}} \dot{u}^\nu ) \mid_{t=0} d\mu_{\tilde{g}} } \vspace{1mm} \\
 & = & -{\displaystyle \int \<(h-\pi_g h)\nu_H, H_{f_g}\> d\mu_g}\vspace{1mm} \\
 & = & -{\displaystyle \int_{f_g(M)} h^2 d\mu_{\tilde{g}} +
\frac{1}{\mu_g} \left( \int_{f_g(M)} h \, d\mu_{\tilde{g}} \right)^2}\, .
\end{array}
$$
By the extreme case of the Cauchy-Schwarz inequality, we conclude that
$h$ is constant. Thus, the function $\| H_{f_g}\|^2$ is constant, and by 
the first identity in (\ref{gr32}), if $n\geq 3$, 
$f_g$ must be a critical point of $\Psi_{f_g}/\mu_g^{2/N}$ in the said space 
of deformations. 

If $M$ is not orientable, we carry the argument above on an orientable 
$2$-to-$1$ cover with the lifted metric, which is locally isometric to the
metric on the base, and so its Nash isometric embedding is locally isometric to
that of $(M,g)$. Since the mean curvature function is constant on the cover, 
$\pm h$ will be a locally defined constant downstairs, and therefore, 
$h^2$ will be a globally well defined constant. The remaining part of the 
argument now follows verbatim the one given above in the orientable case.
\qed
\medskip

{\it Proof of Theorem \ref{th7}}. If $g$ has constant scalar curvature, then 
by (\ref{sce}) the function $\| \alpha_{f_g}\|^2-\| H_{f_g}\|^2$ is constant,
and, by (\ref{gr31}), the embedding $f_g$ is a critical point of 
$\mc{S}_g/\mu_g^{\frac{2}{N}}$ in the space of volume preserving conformal
deformations of it. By Theorem \ref{thgr10}, the function $\| H_{f_{g}}\|^2$ 
is constant, and $f_{g}$ is a critical point of 
$\Psi_{f_g}(M)/{\mu^{\frac{2}{N}}_{g}}$ in the said space of deformations.
Hence, the function $\| \alpha_{f_{g}}\|^2$ is constant also, and by 
(\ref{gr32}), $f_g$ is a critical point of the functional  
$\Pi_{f_g}(M)/{\mu^{\frac{2}{N}}_{g}}$ as well. 

The converse is straightforward. If the functions $\| H_{f_g}\|^2$ and 
$\| \alpha_{f_g}\|^2$ are constants, by (\ref{sce}) we conclude that
$s_g$ is constant, and, by (\ref{lu1}) and (\ref{gr32}) that the isometric 
embedding $f_g$ is a critical point of the Yamabe functional in the space
of volume preserving conformal deformations of it. 
\qed 
\medskip

{\it Proof of Theorem \ref{newth8}}. If $t\rightarrow f_{g_t}$ is the path
of Nash isometric embeddings of the path of metrics $t\rightarrow g_t$, we
can produce a two parameter family of functions $(t,s) \rightarrow v_t(s)$ such
that the metric $g_{s,t}=e^{2v_t(s)}\tg$ has associated family of embeddings
$f_{g_{s,t}}$, which for each fixed $t$, deforms conformally $f_{g_t}$ at $s=0$ 
into $f_{g_t^Y}$ at $s=1$. By (\ref{eq3}), the constant function
$\| H_{f_{g_{s,t}}}\|^2$ depends on
$\| H_{f_{g_t}}\|^2$, $v_t(s)$, and the normal component of 
$\nabla^{\tg}v_t(s)$. Hence, $\| H_{f_{g_{1,t}}}\|^2=\| H_{f_t^Y}\|^2$ is as 
regular, as a function of $t$, as are the metrics $g_t$.

We consider a path of normal vector fields $\nu_{H_{f_t^Y}}$ such that 
$H_{f_{g_t^Y}}=h_{g_t^Y}\nu_{H_{f_{g_t^Y}}}$, with  
$t \rightarrow h_{f_{g_t^Y}}$ a path of constant functions. 
By Theorem \ref{th7}, $\| H_{f_{g_t^Y}}\|^2$ and
$\| \alpha_{f_{g_t^Y}}\|^2$ are both constant functions, and 
$f_{g_t^Y}$ is a critical point of the functionals $\Psi_{f_g}(M)$ 
and $\Pi_{f_g}(M)$ in the space of volume preserving conformal deformations of 
the embedding, respectively. Since homothetics of Yamabe metrics remain
Yamabe metrics, the isometric embedding of $(s,t) \rightarrow e^{2s}g_t^Y$ 
remains a critical point of these two extrinsic functionals in 
the said space of conformal deformations. The embedding 
$f_{e^{2s}g_t^Y}$ becomes a critical point of $\Psi_{f_g}(M)$ in 
the whole space of its conformal deformations if we choose
$s$ to make $\Psi_{f_{e^{2s}g_t^Y}}(M)$ stationary 
along the $H_{f_{g_t^Y}}$ direction also, lifting the restriction on the
volume. 

We let $e^{u_t(s)}g_t^Y$ be a path of conformal dilation deformations of 
$g_t^Y$, that is, $u_t(s)=s$ on points of $f_{g_t^Y}(M) \hookrightarrow
\mb{S}^{\tn}$. By (\ref{eq5}), $s_{e^{2s}g_t^Y}=e^{-2s}s_{g_t^Y}$. Hence, 
by the strong identity (\ref{lu2}), we have that
$$
n=e^{-2s}(n +\tilde{g}(2H_{f_{g_t^Y}}-n\nabla^{\tilde{g}}u_t^\nu, \nabla^{\tg} 
u_t^\nu ))\, , 
$$
and by (\ref{eq3}), we conclude that
$$
\| H_{f_{e^{2s}g^Y}}\|^2 = e^{-2s}( h_{f_{g^Y_t}}^2 - n^2(e^{2s}-1) )\, .
$$
The function 
$$
s\rightarrow \Psi_{f_{e^{2s}g_t^Y}}(M)= \int e^{(n-2)s}( h_{g^Y_t}^2
- n^2(e^{2s}-1) ) d\mu_{g_t^Y}
$$
has two critical points defined by the values of $s$ such 
that $h_{g_t^Y}^2/n=n(e^{2s}-1)$, and $(n-2)(h_{g_t^Y}^2+n^2) = n^3e^{2s}$, 
respectively. The first critical point corresponds to the absolute 
minimizer of $\Psi_{f_{e^{2s}g_{t^Y}}}(M)$, with 
 $e^{2s}= (n^2+h_{g_t^Y}^2)/n^2$, $\nabla^{\tg}u_t(s)^{\nu}\mid_{f_{g_t^Y}(M)}=
H_{f_{g_t^Y}}$, and by (\ref{eq4}), $\| \alpha_{f_{e^{2s}g_t^Y}}\|^2 = 
e^{-2s}\left(\| \alpha_{f_{g_t^Y}}\|^2 -\frac{h_{g_t^Y}^2}{n}\right)$.
The latter critical point has nonzero mean curvature function and  
$e^{2s}\neq 1$ always. 
\qed  
\medskip

Theorem \ref{mt} is a strengthened version of \cite[Theorem 2]{rss2}, based
on the now available Theorem \ref{thgr10}, and the Palais' isotopic 
deformations of the embeddings.
\medskip

{\it Proof of Theorem \ref{mt}.} If we write $H=h\nu_{H}$, since 
$h$ is constant up to a sign, if the first inequality in (\ref{es}) holds 
for any $\lambda \in [0,\frac{1}{2})$, then
$$
\left( \frac{1}{2}-\lambda\right) h^2 \leq n -\| \nabla^\nu_{e_i} \nu_H\|^2 -
\frac{1}{2}h^2 + {\rm trace}\, A_{\nu_H}^2  \, ,
$$
and the left hand side cannot be zero if $h\neq 0$. By the  
 critical point equation \cite[Theorem 3.10]{gracie}, we obtain that
$$
0=2\int h\Delta h d\mu_g = \int h^2( 2n- 2\| \nabla_{e_i}^\nu \nu_H\|^2 -h^2 +
2 {\rm trace}\,A_{\nu_H}^2) d\mu_g  
$$
can only be zero if $h=0$. It follows that the embedding is minimal, and 
satisfies the hypothesis of the
standard gap theorem of Simons \cite[Theorem 5.3.2, Corollary 5.3.2]{si} 
\cite[Main Theorem]{cdck} \cite[Corollary 2]{bla}, which proves (1) in its
entirety. 

If the first inequality in (\ref{es}) fails to hold for any
$\lambda \in [0,\frac{1}{2})$, by the critical point equation 
\cite[Theorem 3.10]{gracie}, we have that
$$
\left( \frac{1}{2}-\lambda\right) h^2 > n -\| \nabla^\nu_{e_i} \nu_H\|^2 -
\frac{1}{2}h^2 + {\rm trace}\, A_{\nu_H}^2 =0 \, ,
$$
so the embedding is not minimal, and we have that
$$      
0\leq {\rm trace}\,A_{\nu_H}^2= \frac{1}{2}h^2+\| \nabla_{e_i}^\nu \nu_H\|^2 -n \leq \| \alpha\|^2_{f_g}\, ,
$$ 
the first and last of these inequalities because $A_{\nu_H}^2$ is nonnegative
with trace bounded above by $\| \alpha_{f_g}\|^2$. 
No $C^2$ Palais deformation of $f_g$ by normal 
stationary critical embeddings, which maintains this nonnegativity 
throughout, can ever be minimal.
\qed
\medskip

{\it Proof of Theorem \ref{nth9}.} Since $f_g: (M,g) \hookrightarrow 
(\mb{S}^{n+1},\tilde{g})$ is minimal, and $s_g=n(n-2)$, by Theorem \ref{th4}, 
$g$ is a Yamabe metric in $[g]$ (see Corollary \ref{coam}). We have that
$2(n-1)W_g(\omega^\sharp_g,\omega_g^\sharp)=n^2$, 
so,  by (\ref{eqn3}), $(n-1)s^J_{g}=2n(n-1)$, and by Theorem 
\ref{t21}, $g$ is not a $J$ Yamabe metric. We choose and fix a $J$
Yamabe metric $g^{JY}$ in $[g]$, of the same volume $\omega_n$ as that of the 
standard sphere $\mb{S}^n\subset \mb{R}^{n+1}$. 

Since $(M,J,g^{JY})$ does not achieve the universal bound of the conformal 
invariant $\lambda^J(M,[g])$, and $g$ is a Yamabe metric, by the first of the 
inequalities in (\ref{10}), we have the estimates
$$
n\frac{n-2}{n-1} \left( \frac{\mu_g(M)}{\omega_n}\right)^{\frac{2}{n}} <
s^J_{g^{JY}} < n \, , 
$$  
where regardless of which of the three $(M,g)$s is under consideration, we 
have that 
$$
0.91 < \frac{n-2}{n-1} \left( \frac{\mu_g(M)}{\omega_n}\right)^{\frac{2}{n}} 
< 0.97 \, . 
$$

If $t \rightarrow f_{\tg_t^Y}$ is assumed to vary in (at least) a $C^2$ manner,
since $\lambda^J(M,[g])$ varies continuously with the conformal class 
(see Lemma 
\ref{l44}), there exists a path of functions $t \rightarrow u_t$ such that
$g_t^{JY}= e^{2u_t}\tilde{g}_t^Y$ is $J_t$ Yamabe of volume $\omega_n$, 
$g_0^{JY}=g^{JY}$. By (\ref{eq6}), the $J_t$ scalar curvature of this path
satisfies the equation
$$
(n-1)s_{g_t^{JY}}^{J_t}=e^{-2u_t}(n(n-1)-\| \alpha_{f_{\tilde{g}_t^Y}} \|^2+
2(n-1)W_{\tilde{g}_t^Y}(\omega_{\tilde{g}_t^Y}^\sharp, 
\omega_{\tilde{g}_t^Y}^\sharp) + (n-1)(2\Delta^{\tilde{g}_t^Y}u_t - 
(n-2)\tilde{g}_t^Y(\nabla^{\tilde{g}_t^Y}u_t^\tau,\nabla^{\tilde{g}_t^Y}
 u_t^\tau))) \, .
$$
By integration with respect to the measure $e^{\frac{n+2}{2}u_t}d\mu_{g_t^Y}$,
we obtain that 
$$
(n-1)s_{g_t^{JY}}^{J_t}\int e^{\frac{n+2}{2}u_t}d\mu_{g_t^Y} =
\int (n(n-1)-\| \alpha_{f_{\tilde{g}_t^Y}} \|^2+
2(n-1)W_{\tilde{g}_t^Y}(\omega_{\tilde{g}_t^Y}^\sharp, 
\omega_{\tilde{g}_t^Y}^\sharp))e^{\frac{n-2}{2}u_t}d\mu_{g_t^Y}\, ,  
$$
and since at $t=0$, $\tilde{g}_t^{Y}=g$, $J_t=J$, and 
$(n(n-1)-\| \alpha_{f_{\tilde{g}_t^Y}} \|^2+ 2(n-1)W_{\tilde{g}_t^Y}(
\omega_{\tilde{g}_t^Y}^\sharp, \omega_{\tilde{g}_t^Y}^\sharp) = n(n-1)$,  
we have that 
$$
\frac{s^J_{g^{JY}}}{n} = \left( {\int e^{\frac{n-2}{2}u_t}d\mu_{g_t^Y}}
\middle/ {\int e^{\frac{n+2}{2}u_t}d\mu_{g_t^Y}}\right) \mid_{t=0} < 1 
$$
By continuity, there exists an $\varepsilon'$,
$0< \varepsilon' < \varepsilon$ such that 
$$
\frac{s^{J_t}_{g^{JY}_t}}{n} = \left( {\int e^{\frac{n-2}{2}u_t}d\mu_{g_t^Y}}
\middle/ {\int e^{\frac{n+2}{2}u_t}d\mu_{g_t^Y}}\right) < 1 
$$
for all $t \in (-\varepsilon', \varepsilon')$. 

On the other hand, by the strong identity (\ref{lu2}), we have that
$$
n(n-1)=e^{-2u_t}(n(n-1)+(n-1)(2\Delta^{\tilde{g}_t^Y}u_t - 
(n-2)\tilde{g}_t^Y(\nabla^{\tilde{g}_t^Y}u_t^\tau,\nabla^{\tilde{g}_t^Y}
 u_t^\tau))-n(n-1) \tilde{g}_t^Y(\nabla^{\tilde{g}_t^Y}u_t^\nu,
\nabla^{\tilde{g}_t^Y} u_t^\nu) \, , 
$$
and, by (\ref{eq3}), that
$$
h_{f_{g_t^{JY}}}^2 = e^{-2u_t}n^2 \tilde{g}_t^Y(\nabla^{\tilde{g}_t^Y}u_t^\nu,
\nabla^{\tilde{g}_t^Y} u_t^\nu) \, . 
$$
By Theorem \ref{thgr10}, this extrinsic function is constant. Hence, 
integrating the previous identity with respect to the measure 
$e^{\frac{n+2}{2}u_t} d\mu_{g_t^Y}$, we obtain that
$$
n(n-1)\int e^{\frac{n+2}{2}u_t} d\mu_{g_t^Y} =   
n(n-1) \int e^{\frac{n-2}{2}u_t} d\mu_{g_t^Y}-\frac{n(n-1)}{n^2}
h_{f_{g_t^{JY}}}^2 \int e^{\frac{n+2}{2}u_t} d\mu_{g_t^Y}\, , 
$$
from which it follows that 
$$
1+\frac{h_{f_{g_t^{JY}}}^2}{n^2} =
 \left( {\int e^{\frac{n-2}{2}u_t}d\mu_{g_t^Y}}
\middle/ {\int e^{\frac{n+2}{2}u_t}d\mu_{g_t^Y}}\right) \geq 1 
$$
for all $t$. 
\qed

\section{Fourth: The second main theorem} 
For no metric in the conformal class of the standard product of 
any of the manifolds $\mb{S}^2(r_1) \times \mb{S}^4(r_2)$, 
$\mb{S}^2(r_1) \times \mb{S}^6(r_2)$, or $\mb{S}^6(r_1) \times \mb{S}^6
(r_2)$, there exists an orthogonal integrable almost complex structure
$J$ \cite{euse}. Indeed, if $J$ were an integrable structure on  
the Riemannian product $\mb{S}^2(r_1) \times \mb{S}^k(r_2)$, $k=2,4,6$,
by the Hermitian identity \cite{gre} 
$$
\sum_{\begin{array}{c}
(i_1,i_2,i_3,i_4)\, , i_k \in \{ 0,1\} \\
\sum i_k = 0,2, 4
\end{array}}
(-1)^{\frac{|i|}{2}}R^{g}(J^{i_1} X, J^{i_2} Y, J^{i_3} Z, J^{i_4} W) = 0 \, ,
$$ 
and the nonnegativity of the curvature in the second factor, 
it would follow that $J$ induces by projection the canonical complex structure 
on the first factor, and an integrable almost complex 
structure on the second. On the other hand, for any orthogonal almost 
complex structure $J$ 
on the product $\mb{S}^6(r_1)\times \mb{S}^6(r_2)$, the $J$-Ricci tensor 
$r_g^J$ is symmetric, $J$ descends to each of the factors, and we have a 
product almost Hermitian structure. 

\begin{theorem}
None of the manifolds $\mb{S}^2 \times \mb{S}^4$, $\mb{S}^2 \times \mb{S}^6$, 
or $\mb{S}^6 \times \mb{S}^6$, respectively, is diffeomorphic to a 
complex manifold.
\end{theorem}

{\it Proof}. We let $(M^n,J,g)$ be any of the Riemannian manifolds 
$(\mb{S}^{6,2},g_{\mb{S}^{6,2}})$, $(\mb{S}^{8,2},g_{\mb{S}^{8,2}})$, or
$(\mb{S}^{12,6},g_{\mb{S}^{12,6}})$, respectively, with their octonionic
induced almost complex structure, and minimal isometric embedding  
$f_g: (M,g) \hookrightarrow (\mb{S}^{n+1},\tilde{g})$. By Theorem \ref{th4},
$g$ is a Yamabe metric in its class $[g]$, of scalar curvature $s_g=n(n-2)$, 
which by Theorem \ref{t21} is not $J$ Yamabe because $2(n-1)W_g(\omega^\sharp_g,
\omega_g^\sharp)=n^2$. We choose, and fix, a $J$ Yamabe metric $g^J$ 
in $[g]$ of the same volume $\omega_n$ as that of the standard sphere 
$\mb{S}^n\subset \mb{R}^{n+1}$. 

\begin{enumerate}
\item[A)] If $J'$ is an almost complex structure on $M^n$ in the same
orientation class as that of $J$, we choose a path 
$[0,1]\ni t\rightarrow J_t$ of almost complex structures on $M$ connecting 
$J'$ and $J$, and consider the path of metrics 
$$
[0,1] \ni t \rightarrow g_t( \, \cdot \, , \, \cdot \, ) = \frac{1}{2}(  
g( \, \cdot \, , \, \cdot \, )+g( J_t \, \cdot \, , J_t \, \cdot \, ))  \, . 
$$
We obtain a path $[0,1]\ni t\rightarrow  (J_t, g_t)$ of 
almost Hermitian structures on $M$ that begins at 
structure $(J',g_0)$, and ends at $(J,g)$, with associated paths
$t\rightarrow [g_t]$ of conformal classes and
$t\rightarrow (J_t, [g_t])$ of almost Hermitian pairs. 

\item[B)] By solving the Yamabe and $J$ Yamabe problems on $[g_t]$ and 
$(J_t, [g_t])$, respectively, we produce a Yamabe metric $g_t^Y$ of scalar 
curvature $s_{g_t^Y}$, and $J_t$ Yamabe metric $g^{JY}_t$ of $J_t$ scalar 
curvature $s_{g_t^{JY}}^{J_t}$, which we normalize to have volume 
$\mu_g(M)$, the volume of $M$ in the metric $g$, and $\omega_n$, respectively.
Since $g$ is a Yamabe metric, 
and $g^J$ is a $J$ Yamabe metric, by Lemma \ref{l44}, we can make the choices 
of $g_t^Y$ and $g_t^{JY}$ to produce (at least) $C^2$ paths of metrics 
$t \rightarrow g_t^Y$, and $t \rightarrow g_t^{JY}$, that equal $g$ and $g^J$ 
at $t=1$, respectively.

\item[C)] By applying the Nash isometric embedding theorem \cite{nash}, for a
sufficiently large $p$, we may obtain a family of isometric embeddings  
$$
f_{g^{Y}_t}: (M,g_t^Y) \hookrightarrow  (\mb{S}^{\tn= n +p },\tilde{g}) \, , 
$$
and
$$
f_{g^{JY}_t} : (M,J_t,g_t^{JY}) \hookrightarrow  (\mb{S}^{\tn= n +p},
\tilde{g})\, ,
$$
such that, in terms of the associated conformal invariants,  
$$
s_{g_t^Y}=\frac{1}{\mu_g^{\frac{2}{n}}}\lambda(M,[g_t])=
n(n-1)-(\| \alpha_{f_{g^{Y}_t}} \|^2- \| H_{f_{g^{Y}_t}}\|^2) \, , 
$$
and 
$$
s_{g_t^{JY}}^J= 
\frac{1}{\omega_n^{\frac{2}{n}}} \lambda^{J_t}(M,[g_t]) 
 =n- 
\< \alpha^{JY}_{f_{g^{JY}_t}} (e_i,J_te_j),\alpha^{JY}_{f_{g^{JY}_t}} 
(e_j,J_te_i)\> \, , 
$$
respectively. By Theorem \ref{thgr10}, the function $\|H_{f_{g^{Y}_t}}\|^2$
is constant.

\item[D)] Up to an isometry of the background sphere, the embedding 
$f_{g^{Y}_t}$ at $t=1$ 
coincides with the standard isometric embedding $f_g: (M,g) \hookrightarrow 
(\mb{S}^{n+p},\tilde{g})$ as a minimal submanifold of scalar curvature 
$n(n-2)$, and we have 
$$
s_g=s_{g_1^Y}=\frac{1}{\mu_g^{\frac{2}{n}}} \lambda(M,[g])=n(n-2) \, .
$$
The embedding $f_{g^{JY}_1}$ is an isometric embedding of $g^J$, 
a metric in $[g]$ other than $g$, and by the positivity of $W_g(\omega_g^\sharp,
\omega_g^\sharp)$, we have that $(n-1)\lambda^J(M,[g]) > \lambda(M,[g])$.  

\item[E)] Since $t\rightarrow [g_t] \rightarrow g_t^Y$ is at (least) $C^2$, 
and $s_{g_t^Y}$ is constant, by (\ref{sce}) the function 
$\| \alpha_{f_{g^{Y}_t}}\|^2-\| H_{f_{g^{Y}_t}}\|^2$ is constant, and the path
$t \rightarrow \| \alpha_{f_{g^{Y}_t}}\|^2-\| H^Y_{f_{g^{Y}_t}}\|^2$ is 
(at least) continuous. 
By Theorem \ref{th7}, the functions $\| \alpha_{f_{g^{Y}_t}}\|^2$
 and $\| H_{f_{g^{Y}_t}}\|^2$ are individually constants, and the embedding 
$f_{g_t^Y}$ is a stationary point of the extrinsic 
functionals $\Psi_{f_g}(M)/\mu_g^{2/N}$, and $\Pi_{f_g}(M)/\mu_g^{2/N}$
in the space of its volume preserving conformal deformations, respectively.  

\item[F)] We write $H_{f_{g_t^Y}}$ as $H_{f_{g^{Y}_t}}=h_t
\nu_{f_{g_t^Y}}$, where $\nu_{f_{g_t^Y}}$ is a normal section of the normal
bundle of $f_{g_t^Y}(M)$ in $\mb{S}^{n+p}$.
 Since the $f_{g^{Y}_t}$s are volume
preserving, the component $T^{\nu}$ of the variational 
vector field $T=T^{\tau}+T^{\nu}$ is $L^2$-orthogonal to
$H_{f_{g^{Y}_t}}$. But normal directions are all conformal, so by (E) above, 
$\Psi_{f_{g_t^Y}}$ and $\Pi_{f^Y_{g_t}}$ are stationary in any
normal direction $L^2_{g_t^Y}$ orthogonal to $\nu_{f_{g_t^Y}}$ also. 
By Theorem \ref{newth8}, the conformally dilated metrics 
$\tilde{g}_t^Y=(1+h_t^2/n^2)g^Y_t$ have their corresponding family of 
isometric embeddings
$$
f_{\tilde{g}_t^Y}: (M,(1+h_t^2/n^2)g_t^Y) \rightarrow (\mb{S}^{n+p},\tilde{g})
$$
minimal for all $t$, so they all are critical points of the functional
$$
f_g \rightarrow  \Psi_{f_g}(M) \, ,
$$
along all normal directions, including $\nu_{f_{g_t^Y}}$ \cite{gracie}, and
$t\rightarrow \| \alpha_{f_{\tg_t^Y}}\|^2$ is continuous.
By the minimality of $f_{g_1^Y}$, we have that $\tilde{g}_1^Y=g^Y_1=g$, and
the embeddings $f_{\tilde{g}_1^Y}$ and $f_{g_1^Y}$ coincide.

\item[G)] At $t=1$, we have that $\| \alpha_{f_{\tilde{g}_1^Y}}\|^2=n$,  
$\| H_{f_{\tilde{g}_1^Y}}\|^2=0$, and the estimates in (\ref{es}) hold. 
By Theorem \ref{mt}, we conclude that, modulo an isometry, $f_{\tilde{g}_1^Y}=
f_{g_1^Y}$ is the standard minimal embedding $f_g : (M,g)
\rightarrow (\mb{S}^{n+p},\tilde{g})$, and we have that $[\tilde{g}_1^Y]=
[g^Y_1]=[g_1]=[g]$, and that this remains true for any $t$ for which 
$\| \alpha_{f_{\tilde{g}_t^Y}}\|^2=n$, in which case, 
$[\tilde{g}_t^Y]=[g^Y_t]=[g_t]=[g]$. Thus, there exists a smallest 
$a$, $0\leq  a \leq 1$ such that if $t\in [a,1]$, $\tilde{g}_t^Y=g_t^Y$,
and $f_{\tilde{g}_t^Y} = f_{g_t^Y}=f_{g}$.
 
\item[H)] There there are no 
orthogonal complex structures on the standard
$(M,g)$ \cite{euse}. Hence, if we assume that 
the path $J_t$
starts at an integrable $J'=J_0$, we must have that
$[{\tilde{g}^Y_0}]=[g^Y_0]=[g_0]\neq [g]$, and by Theorem \ref{mt},
$[g_0] \neq [g]$, and the constant function 
$\| \alpha_{f_{\tilde{g}^Y_0}}\|^2$ must satisfy the inequality
$$
\| \alpha_{f_{\tilde{g}^Y_0}}\|^2 > n \geq \frac{np}{2p-1} \, .
$$
We conclude that the constant $a$ in (G) above is strictly positive, and by 
Theorem \ref{mt} again, that there exists $\varepsilon$, $0< \varepsilon < a$, 
such that for any $t \in (a-\varepsilon, a]$, the scalar curvature 
$s_{\tilde{g}_t^Y}$, the function $W_{\tilde{g}_t^Y}(
\omega_{\tilde{g}_t^Y}^\sharp, \omega_{\tilde{g}_t^Y}^\sharp)$, and the
$J_t$ scalar curvature $s_{g_t^{JY}}^{J_t}$ are 
positive, $[g_t] \neq [g]$ for any 
$t \in (a-\varepsilon, a)$, and in this latter range of $t$s, 
$$
\| \alpha_{f_{\tilde{g}^Y_t}}\|^2 > n \geq \frac{np}{2p-1} \, . 
$$
By Theorem \ref{nth9}, the path $[0,1]\ni t \rightarrow \| 
\alpha_{f_{\tilde{g}^Y_t}}\|^2$ changes 
discontinuously across $t=a$, a contradiction. 
\end{enumerate}
\qed

\end{document}